\documentclass{crm-article}

\usepackage[T1,T2A]{fontenc}
\usepackage[utf8]{inputenc}

\usepackage{nccmath}
\usepackage{mathtext}
\usepackage{amsmath,amssymb}

\usepackage{float}
\usepackage{caption}
\usepackage[export]{adjustbox}
\usepackage{graphicx}

\usepackage{tabularx}
\usepackage{multirow}
\usepackage{multicol}

\usepackage{ragged2e}
\usepackage{pdfpages}

\usepackage[usenames,dvipsnames]{xcolor}
\newcommand{\ok}[2]{\colorbox{ForestGreen!#1}{\color{black}{#2}}}
\newcommand{\indif}[2]{\colorbox{Gray!#1}{\color{black}{#2}}}

\begin{document}

\journalVol{10}
\journalNo{1} 
\setcounter{page}{1}

\journalSection{Математические основы и численные методы моделирования}
\journalSectionEn{Mathematical modeling and numerical simulation}

\journalReceived{??.12.2021.}
\journalAccepted{??.12.2021.}

\UDC{519.856}
\title{Стохастическая оптимизация в задаче цифрового предыскажения сигнала}
\titleeng{Stochastic optimization in digital pre-distortion of the signal}
\thanks{Работа Д.А. Пасечнюка и А.М. Райгородского поддержана грантом Российского научного фонда (проект №21-71-30005).}
\thankseng{The research of D. Pasechnyuk and A. Raigorodskii was supported by Russian Science Foundation grant (project No. 21-71-30005).}

\author{\firstname{А.\,В.}~\surname{Алпатов}}
\authorfull{Артём Вадимович Алпатов}
\authoreng{\firstname{A.\,V.}~\surname{Alpatov}}
\authorfulleng{Artem V. Alpatov}
\email{artemalpatov2@gmail.com}
\affiliation{ГАОУ СО ``Самарский лицей информационных технологий (Базовая школа РАН)'',\protect\\ Россия, 443096, г. Самара, ул. Больничная, д. 14а}
\affiliationeng{``Samara Lyceum of information technologies (Basic school of the RAS)'', Samara, 443096, Russia}

\author[2]{\firstname{Е.\,А.}~\surname{Петерс}}
\authorfull{Егор Александрович Петерс}
\authoreng{\firstname{E.\,A.}~\surname{Peters}}
\authorfulleng{Egor A. Peters}
\email{cool47.cool@gmail.com}
\affiliation[2]{ЦДНИТТ ``УникУм'' при КузГТУ,\protect\\ Россия, 650000, г. Кемерово, ул. Красноармейская, д. 117}
\affiliationeng{CCSETC at KuzSTU ``UnicUm'', Kemerovo, 650000, Russia}

\author[3,4]{\firstname{Д.\,А.}~\surname{Пасечнюк}}
\authorfull{Дмитрий Аркадьевич Пасечнюк}
\authoreng{\firstname{D.\,A.}~\surname{Pasechnyuk}}
\authorfulleng{Dmitry A. Pasechnyuk}
\email{pasechnyuk2004@gmail.com}
\affiliation[3]{Московский физико-технический институт,\protect\\ Россия, 141701,
г. Долгопрудный, 
Институтский пер., д. 9}
\affiliationeng{Moscow Institute of Physics and Technology, Moscow, 141701, Russia}

\affiliation[4]{Исследовательский центр доверенного искусственного интеллекта ИСП РАН,\protect\\ Россия, 109004, г. Москва, ул. А. Солженицына, д. 25}
\affiliationeng[4]{ISP RAS Research Center for Trusted Artificial Intelligence, Moscow, 109004, Russia}

\author[3]{\firstname{А.\,М.}~\surname{Райгородский}}
\authorfull{Андрей Михайлович Райгородский}
\authoreng{\firstname{A.\,M.}~\surname{Raigorodskii}}
\authorfulleng{Andrei A. Raigorodskii}
\email{mraigor@yandex.ru}

\begin{abstract}
В данной статье осуществляется сравнение эффективности некоторых современных методов и практик стохастической оптимизации применительно к задаче цифрового предыскажения сигнала (DPD), которое является важной составляющей процесса обработки сигнала на базовых станциях, обеспечивающих беспроводную связь. В частности, рассматривается два круга вопросов о возможностях применения стохастических методов для обучения моделей класса Винера--Гаммерштейна в рамках подхода минимизации эмпирического риска: касательно улучшения глубины и скорости сходимости данного метода оптимизации, и относительно близости самой постановки задачи (выбранной модели симуляции) к наблюдаемому в действительности поведению устройства. Так, в первой части этого исследования, внимание будет сосредоточено на вопросе о нахождении наиболее эффективного метода оптимизации и дополнительных к нему модификаций. Во второй части предлагается новая квази-онлайн постановка задачи и соответственно среда для тестирования эффективности методов, благодаря которым результаты численного моделирования удаётся привести в соответствие с поведением реального прототипа устройства DPD. В рамках этой новой постановки далее осуществляется повторное тестирование некоторых избранных практик, более подробно рассмотренных в первой части исследования, и также обнаруживаются и подчёркиваются преимущества нового лидирующего метода оптимизации, оказывающегося теперь также наиболее эффективным и в практических тестах. Для конкретной рассмотренной модели максимально достигнутое улучшение глубины сходимости составило 7\% в стандартном режиме и 5\% в онлайн постановке (притом что метрика сама по себе имеет логарифмическую шкалу). Также благодаря дополнительным техникам оказывается возможным сократить время обучения модели DPD вдвое, сохранив улучшение глубины сходимости на 3\% и 6\% для стандартного и онлайн-режимов соответственно. Все сравнения производятся с методом оптимизации Adam, который был отмечен как лучший стохастический метод для задачи DPD из рассматриваемых в предшествующей работе \cite{Pasechnyuk}, и с методом оптимизации Adamax, который оказывается наиболее эффективным в предлагаемом онлайн-режиме.
\end{abstract}

\keyword{цифровое предыскажение}
\keyword{обработка сигнала}
\keyword{стохастическая оптимизация}
\keyword{онлайн обучение}

\begin{abstracteng}
In this paper, we test the performance of some modern stochastic optimization methods and practices in application to digital pre-distortion problem, that is a valuable part of processing signal on base stations providing wireless communication. In first part of our study, we focus on search of the best performing method and its proper modifications. In the second part, we proposed the new, quasi-online, testing framework that allows us to fit our modelling results with the behaviour of real-life DPD prototype, retested some selected of practices considered in previous section and approved the advantages of the method occured to be the best in real-life conditions. For the used model, maximum achieved improvement in depth was 7\% in standard regime and 5\% in online one (metric itself is of logarithmic scale). We also achieved a halving of the working time preserving 3\% and 6\% improvement in depth for the standard and online regime, correspondingly. All comparisons are made to the Adam method, which was highlighted as the best stochastic method for DPD problem in paper \cite{Pasechnyuk}, and to the Adamax method, that is the best in the proposed online regime.
\end{abstracteng}

\keywordeng{digital pre-distortion}
\keywordeng{signal processing}
\keywordeng{stochastic optimization}
\keywordeng{online learning}

\maketitle

\begin{figure}[H]
    \centering
    \includegraphics[width=0.3\linewidth]{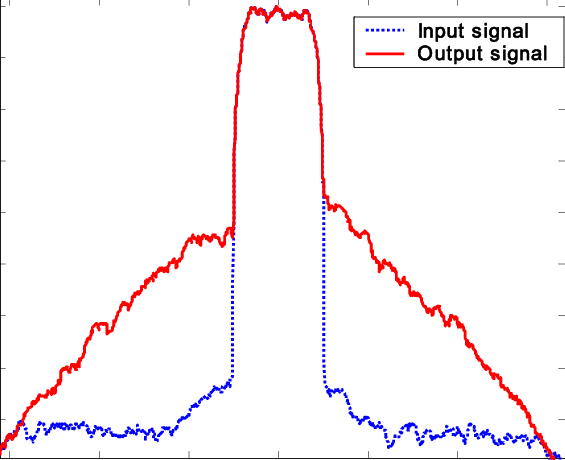}
    \caption{
    Power spectral density plot of original signal (Input signal, blue dotted line) and power amplified signal (Output signal, red line).
        Picture from Wisell~D.\,H. Exploring the sampling rate requirements for behavioural amplifier modelling~// XVIII IMEKO World Congress.~--- 2006.~--- p. 1-4
    }
    \label{fig:pa-harmonics}
\end{figure}

\paragraph{Introduction}
Consider a base station providing some form of wireless communication. Let it be required to transmit with its help some given signal. The signal is represented as a sampled array $x \in \mathbb{C}^m$. One of the transformations that must be applied to the signal before transmission is amplification. To carry out this transformation, the PA device is most often used, whose action on the signal, ideally, has the form $\mathcal{PA}(x) = a \cdot x$, where $a \gg 1$. However, hardware imperfections lead to a deviation from this law; in particular, the $\mathcal{PA}(x)$ signal has out-of-band spurious harmonics, resulting in interference. This is illustrated in Fig.~\ref{fig:pa-harmonics}.

\begin{figure}[H]
    \centering
    \includegraphics[width=0.6\linewidth]{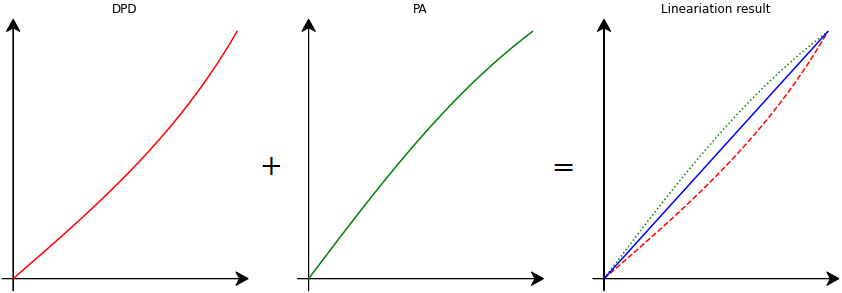}
    \caption{
    The intuitive principle of digital pre-distortion: DPD acts as an inverse of power amplifiers function, so that the resulting transformation is near linear (identical)
        Picture from \cite{Pasechnyuk}
    }
    \label{fig:pa-invert-dpd}
\end{figure}

To prevent noise in the transmitted signal, there is digital pre-distortion technology. Specifically, there is a DPD device acting on the signal just before the PA action, so that the result is the signal $\mathcal{PA}(\mathcal{DPD}(x))$. If we omit the factor $a$, which is easy to accomplish by scaling the signal, then we can see that for an ideal action it must hold that $\mathcal{DPD} = \mathcal{PA}^{-1}$ (Fig.~\ref{fig:pa-invert-dpd}). Thus, the task is reduced to the inversion of $\mathcal{PA}$ function. However, the $\mathcal{PA}$ function actually depends on the characteristics of the environment, the signal x itself, and on time, so it can only be inverted numerically. To do this, we choose a certain parametric family of functions $\{\mathcal{DPD}_\theta\}_{\theta \in \Theta}$ to describe DPD and choose such a parametrization $\theta$ in which $\mathcal{PA} \circ \mathcal{DPD}_\theta$ is  closest to the identical transformation.

\begin{figure}[H]
    \centering
    \includegraphics[width=0.6\linewidth]{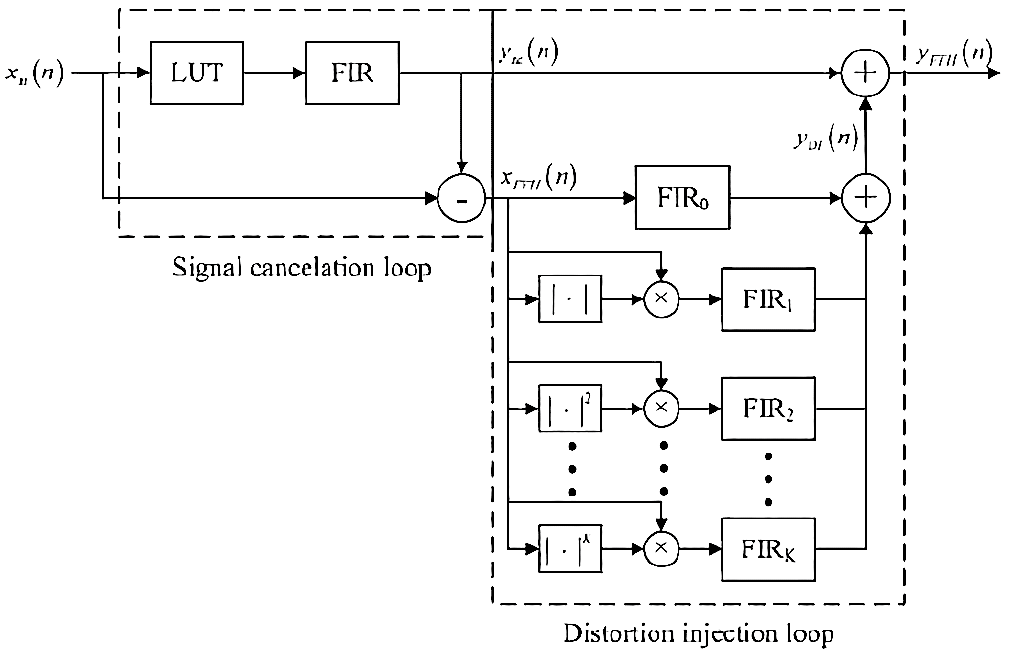}
    \caption{
    The structure of the Wiener--Hammerstein model cell with polynomial non-linearity.
        Picture from \cite{Ghannouchi}
    }
    \label{fig:wiener-hammestein}
\end{figure}

In this study, to model the behaviour of DPD we use the Wiener--Hammerstein cascade models family \cite{Ghannouchi}. The structure of these models is very similar to the structure of multilayer neural network, but instead of the neurons it consists of cells specific for signal processing. Fig.~\ref{fig:wiener-hammestein} shows the principle scheme of such a cell with polynomial non-linearity within. 

Now, when we defined the parametrization for DPD, we should provide a more optimization-friendly formulation to the problem of inverting the $\mathcal{PA}$ function. The point is that to do this we need to have the perfect pre-distortion $\overline{y}$ for some signal $x$, to train the model to repeat this predistortion, i.e. exactly invert the $\mathcal{PA}$ function. Using the standard empirical risk minimization framework, we then obtain some optimization problem with the sum-like function. The following formulation, proposed and described in details in \cite{Pasechnyuk}, is the central one for this paper (we denote the $k$-th element in array $x$ by $[x]_k$ and $x_k$).
\begin{equation} \label{eq:problem}
   \frac{1}{m} \sum_{k=1}^m |[\mathcal{DPD}_\theta(x)]_k-\bar{y}_k|^2\,\rightarrow\,\min_\theta.
\end{equation}

So, the both two research questions of this study correspond to this problem formulation and concern how to solve it properly, taking into account its original motivation.
\begin{itemize}
    \item[{$\boldsymbol{RQ_1}$}] Which stochastic optimization method is the most suitable to solve the stated problem with such a specific data and model behind?
\end{itemize}
In broad strokes, this question was considered in course of study \cite{Pasechnyuk} for the different groups of numerical methods. Here we focus on the stochastic methods (indeed, \eqref{eq:problem} can be easily randomized by terms), and consider several recently proposed approaches in this field.
\begin{itemize}
    \item[{$\boldsymbol{RQ_2}$}] How one should organize the training procedure and evaluation of the model to agree with the practical behaviour of DPD device?
\end{itemize}
The operation of the full DPD technology stack is more complicated than it is representated in \eqref{eq:problem}, and the main difference is its online nature. It is even more online than it is implied, for example, in the expectation minimization problems, because the source of signal (and hence a joint distribution of $(x, \overline{y})$, continuing the analogy) also changes in time. We propose a way to simulate this specificity, and eliminate contradictions between our experimental results and some results of laboratory tests with DPD device.


\paragraph{Standard learning framework}

In this section we consider the standard setup for the problem stated above, i.e. train the model on some known data (it is important for signal processing task, that this subset of the whole data is not random chosen but is the prefix of the ordered signal array) and then test the model on the whole signal (NB: not the remaining part, as usual). In our experiments, length of training signal prefix is 75\% of the whole signal's length (albeit it is deliberately more than needed, see \cite{Pasechnyuk} for details, to avoid any approximation issues that are out of the scope this time). We use the segment with $m=2 \cdot 10^5$ ticks from the real-life 80 MHz signal as a dataset. As we shall see, this way to train/test model is slightly naïve to represent the behaviour of the real-live DPD device (being standard and widely used nonetheless), but thanks to it we can easily test many possible techniques in the simplest (and computationally inexpensive) conditions to select the few best to work with further.

To assess the results of evaluating model on test split of dataset we use the following logarithmic scaling analogue of Mean Squared Error metric with normalization, that is specific to the signal processing tasks:
\begin{equation*}
   \text{NMSE}(y, \overline{y}):=10 \cdot \log_{10} \frac{\sum_{k=1}^{m}|y_k-\overline{y}_k|^2}{\sum_{k=1}^{m}|x_k|^2} \text{ dB}.
\end{equation*}

Source code of the model used to obtain the following experiments results is available at \url{https://github.com/dmivilensky/Sirius-2021-DPD-optimization} (Python 3 + PyTorch). All the computations were performed on the standard personal computer without GPU acceleration.

\subparagraph{Comparison of algorithms}
\begin{figure}[H]
    \centering
    \includegraphics[width=0.6\linewidth]{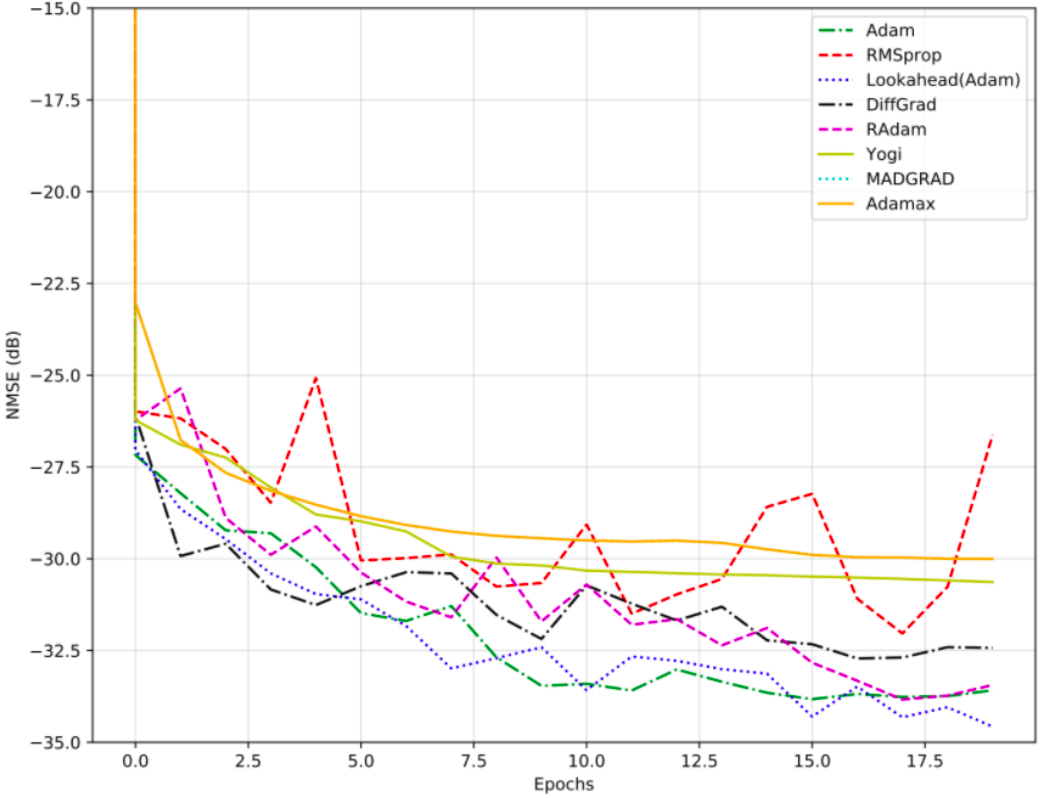}
    \caption{
    The convergence curves of different optimization methods in standard regime at initial epochs
    }
    \label{fig:all-methods}
\end{figure}

\begin{figure}[H]
    \centering
    \includegraphics[width=0.5\linewidth]{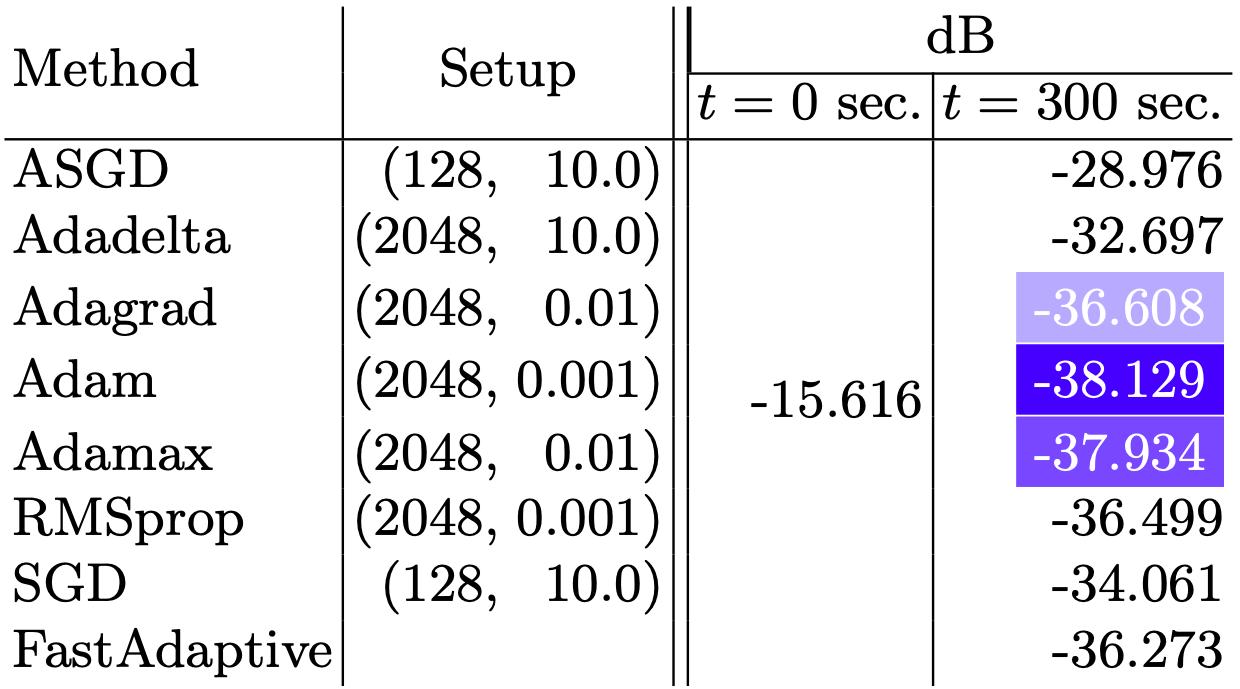}
    \caption{
    Summary of the result of training the DPD model with different optimization methods in standard setting.
        Table from \cite{Pasechnyuk}
    }
    \label{fig:all-methods-table}
\end{figure}

Let us start with the testing of different optimization in application to this standardly formulated learning problem. The convergence of some most known ones was examined in \cite{Pasechnyuk} recently, so there is no need for us to present the same results for our model: in comparison with table, that is shown in Fig.~\ref{fig:all-methods-table}, they are just transposed a couple of decibels higher (in our research, we use the Wiener--Hammerstein model with significantly less parameters than for the model considered in \cite{Pasechnyuk}, so our score will be a little worse~--- it is just a specific of particular used model).

Besides, in this section we focus on the starter methods, i.e. that ones the convergence of which is the most efficient among all tested methods at the beginning of the training process. That's why we show the convergence of all the methods on the first 20 or 40 epochs. The first reason to this is that long term convergence in standard training regime is already considered in \cite{Pasechnyuk}, so we should clarify the short term one, that was not yet considered (it is important for the multi-method optimization schemes). The second is that the convergence on the late epochs and the depth of finally obtained minimum is more the question of late convergence rate of method, not convergence depth, and we will see in online learning framework, that methods presented as the best here are not well adopted to converge with a good rate at late epoch, that is important in practice, because the optimal point we seek for is changing together with data segment. So, let us avoid all this complicated questions in this section, and focus on the starting convergence, i.e. the behaviour of methods far from local minimum.

\begin{figure}[H]
    \centering
    \includegraphics[width=0.5\linewidth]{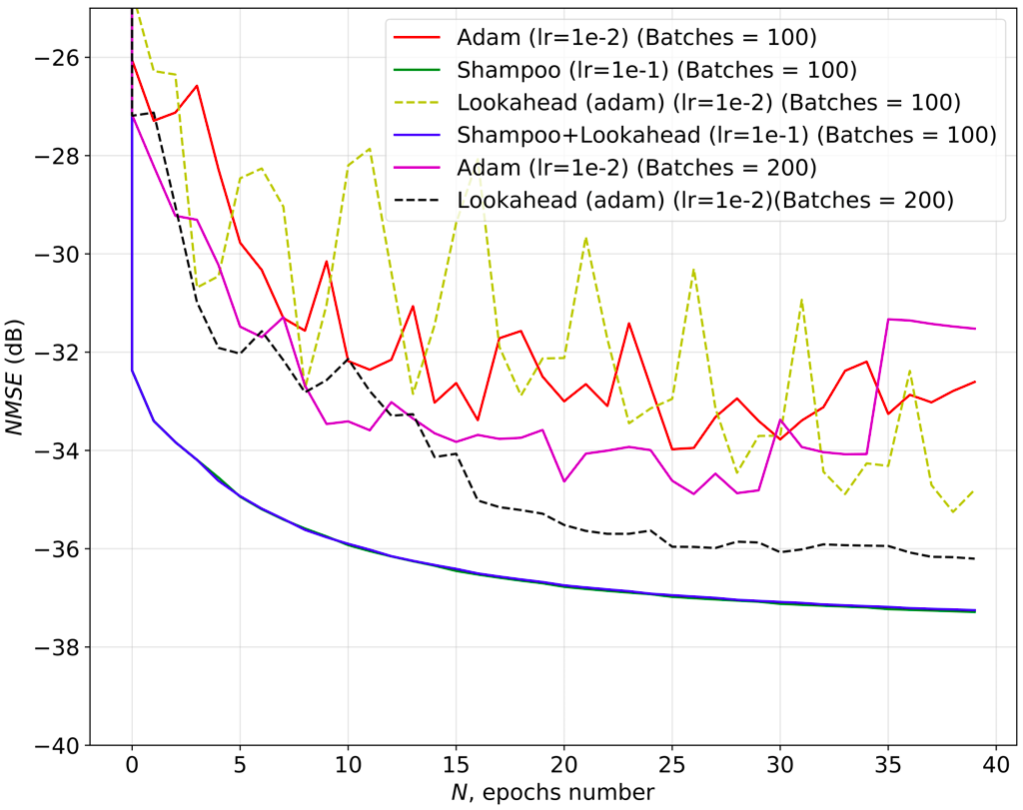}
    \caption{
    The comparison of Adam and Shampoo \cite{Gupta} optimization methods, and the effect of Lookahead envelope on training the DPD model
    }
    \label{fig:lookahead}
\end{figure}

\begin{figure}[H]
    \centering
    \includegraphics[width=0.5\linewidth]{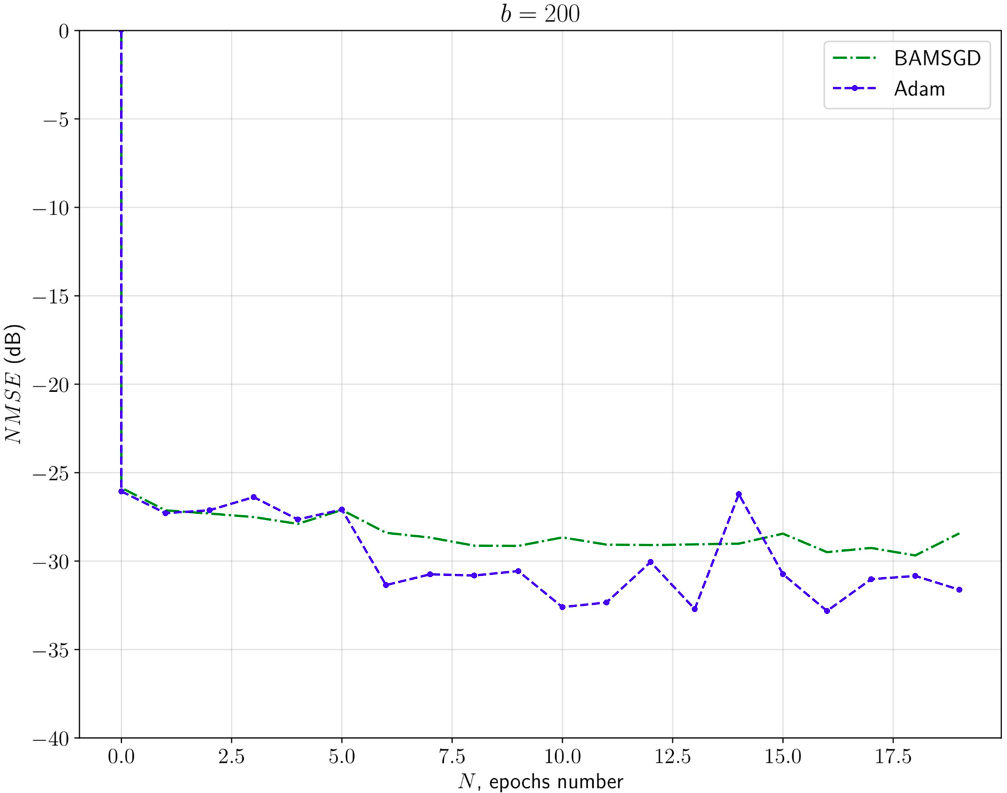}
    \caption{
    The comparison of Adam and AccMbSGD \cite{Woodworth} optimization methods in application to training the DPD model
    }
    \label{fig:accmbsgd}
\end{figure}

So, in Fig.~\ref{fig:all-methods}, Fig.~\ref{fig:lookahead} and Fig.~\ref{fig:accmbsgd} we present the convergence curves of different tested optimization methods. Further, we describe our observations on it.

\begin{enumerate}
    \item Starting convergence of Adamax is not very well, despite the fact that it is one of the most efficient methods in long term perspective. Some of the methods, like RMSprop, are just unstable for the model and signal of used size.
    
    \item The best among methods from Fig.~\ref{fig:all-methods} are Adam, as expected, and Lookahead(Adam) (we also denote it by LaAdam), i.e. the Adam enveloped with Lookahead algorithm from \cite{Zhang}. The latter demonstrates the faster convergence.
    
    \item The only method that outperforms Adam is Shampoo \cite{Gupta}, it is the best starter method for DPD training for now.
    
    \item Lookahead gives very notable results. It significantly improves and stabilizes the convergence of Adam, that is shown in Fig.~\ref{fig:lookahead}. Unfortunately, there is no effect of Lookahead on Shampoo. Nevertheless, convergence of Shampoo is phenomenally stable without Lookahead too. 
\end{enumerate}

Another novel optimization method we tested is Accelerated mini-batch SGD version from \cite{Woodworth} (we name it AccMbSGD). This method is specialized for the overparametrized optimization problems. The point is that in our case the number of parameters in model is significantly bigger than number of terms in sum in empirical risk, so it is reasonable to assume that value of NMSE in minimum is near the zero (and it is). This is the characteristic property of overparametrization, and mentioned method is theoretically optimal for this case.

\begin{enumerate}
    \setcounter{enumi}{4}
    
    \item Unfortunately, the convergence of AccMbSGD is not so well as expected. Moreover, this method is complicated in hyperparameter tuning. It seems that it is not adopted for DPD training. The corresponding curves are presented in Fig.~\ref{fig:accmbsgd}.
\end{enumerate}

\subparagraph{Dynamic batch size}
It is well-known practice to use mini-batching in stochastic optimization methods to improve the convergence (due to the less variance of oracle) and speed up the calculations (it is typical for training neural networks). Usually, the size of the batch is just a hyperparameter of the method and is the same at the every iteration. Nevertheless, it is reasonable to strive to increase of it, if the latter is computationally efficient. On the other hand, this heuristic is quite legitimate: indeed, at the late steps of the method current point is already close enough to minimum and the step direction efficiency becomes more sensitive to the inaccuracies of the stochastic gradient approximation (which just can be reduced by increasing the size of the batch). 

Some theoretical justifications on this technique are provided in \cite{Zhao}. The practical efficiency of batch size adaptation is also shown in \cite{Devarakonda} for some deep learning problems. So, further we check if the changing of batch size affects the convergence of Adam optimization method in our case. We considered both increasing and decreasing changing, as well as the different changing rates.

\begin{figure}[H]
    \centering
    \subfloat[The dependency of NMSE value (dB) on the number of training epochs]{
        \label{fig:dynamic_batch_epochs}
        \includegraphics[width=0.47\textwidth]{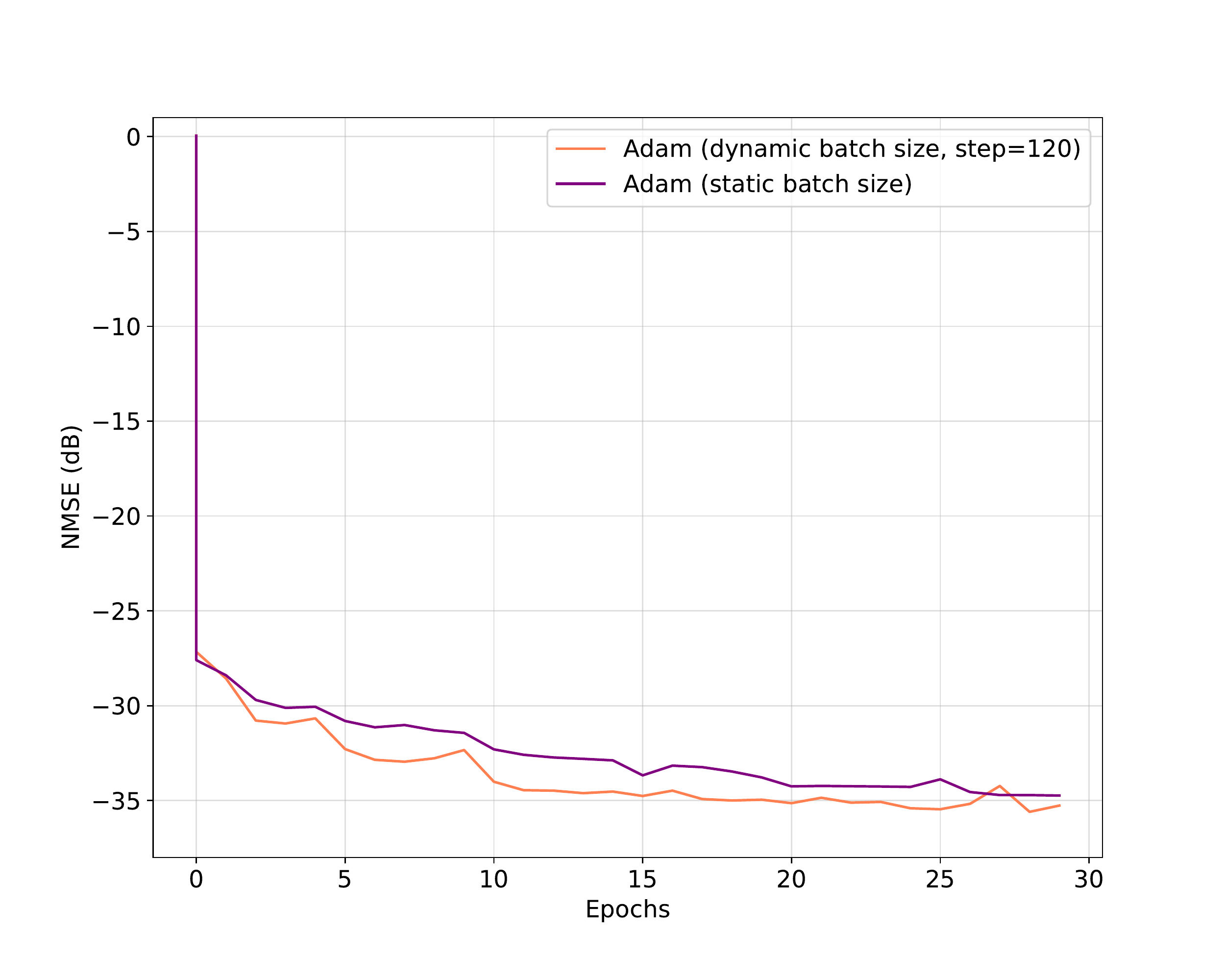}
    }
    \hspace{0.02\textwidth}
    \subfloat[The dependency of working time (s) on the number of training epoch]{
        \label{fig:dynamic_batch_time}
        \includegraphics[width=0.47\textwidth]{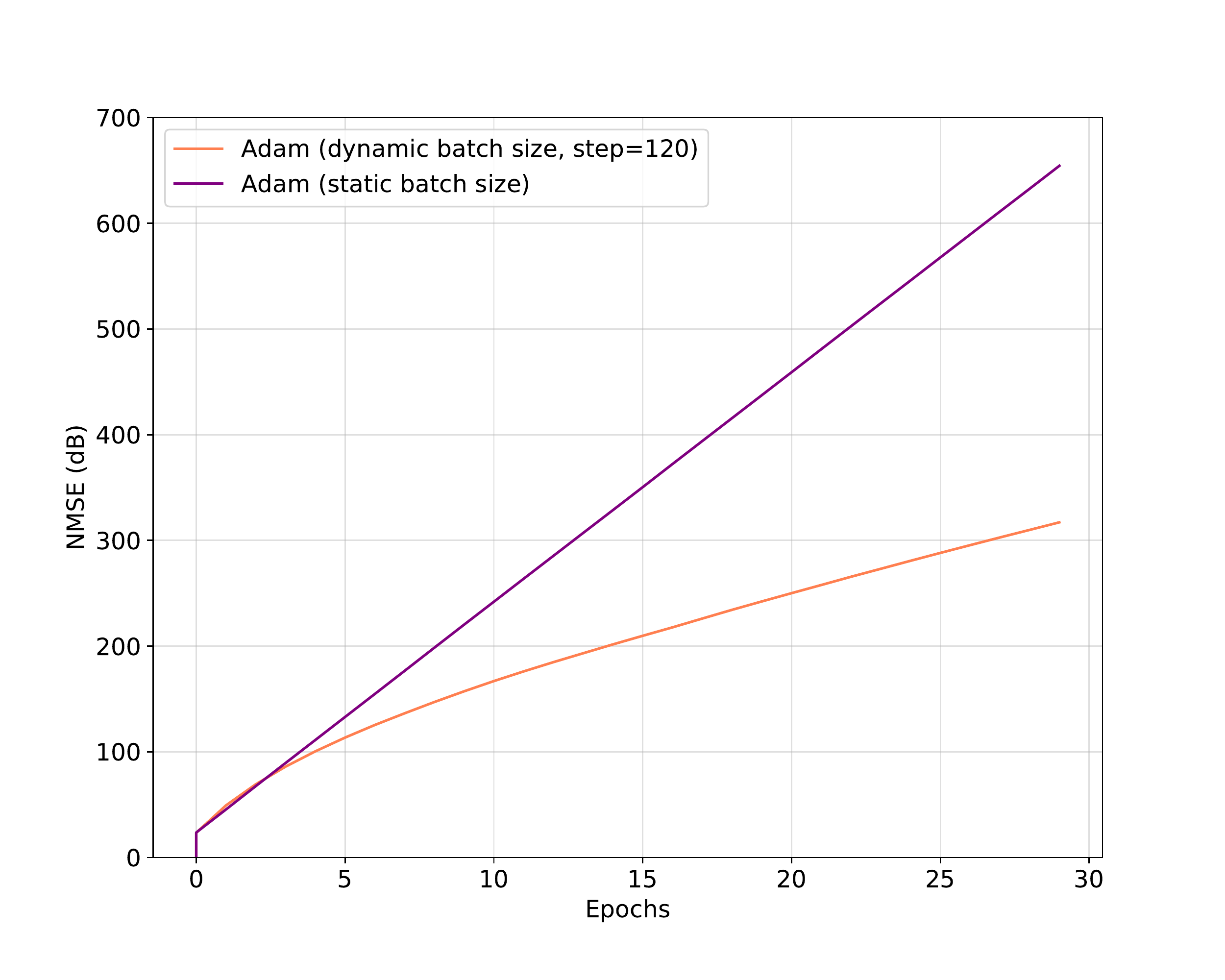}
    }
    \caption{The effect of batch size changing on training the DPD model. It leads to the more deep convergence of Adam (left plot) and slower growth in working time (right plot)}
\end{figure}

\begin{enumerate}
    \item In our experiments, Adam with reducing batch size demonstrated less stable and less deep convergence than the Adam with fixed batch size (as expected).
    
    \item The growth of batch size following the exponential law lead to the improvement of convergence, but its effect is less than for linear law.
    
    \item The best formula for changing the batch size in our case is $b = 200 + 120 \cdot \text{epoch}$. Using it, we obtained a 3\% improvement in NMSE (Fig.~\ref{fig:dynamic_batch_epochs}) and reduced training time by 2 times (Fig.~\ref{fig:dynamic_batch_time}). The results are summarized in Table~\ref{tab:summary_dynamic_batch}.
    
    \item The positive effect in training time increases with an increase in the coefficient in the latter linear formula (it is $120$ there). But here we faced with the trade-off, because the increase in the coefficient can also worsen the limit depth of convergence (value of NMSE). In signal processing task the learning time is minor metric, so it is reasonable to stick around $120$.
\end{enumerate}

\begin{table}[H]
    \centering
    \begin{tabular}{l|c|c}
             & NMSE (dB) & time (s) \\ \hline
    fixed    & -34.90    & 660 \\
    changing & \ok{40}{-35.92} & \ok{40}{320}
    \end{tabular}
    
    \caption{Summary of the result of training the DPD model with fixed and with changing batch size}
    \label{tab:summary_dynamic_batch}
\end{table}

In the experiments above we used the Adam optimization method to demonstrate the efficiency of dynamic batch size approach. The Adam was chosen as the leader among the methods for standard learning framework \cite{Pasechnyuk}. Nevertheless, these results are similar for all the tested optimization methods (they are listed in Table~\ref{tab:methods}).

\subparagraph{Dynamic learning rate}
Another important hyperparameter of the optimization methods is its learning rate, i.e. some starting value of the step size. Further, we describe the results of our experiments on scheduling (changing) the value of learning rate during the training. It is known that for some specific models it may be essential to use a proper learning rate scheduler to obtain the satisfactory accuracy \cite{Xiong}. So, we tested some of the widely used schedulers in our case to check if it affects convergence significantly.

\begin{figure}[H]
    \centering
    \subfloat[The dependency of NMSE value (dB) on the number of training epochs]{
        \includegraphics[width=0.47\textwidth]{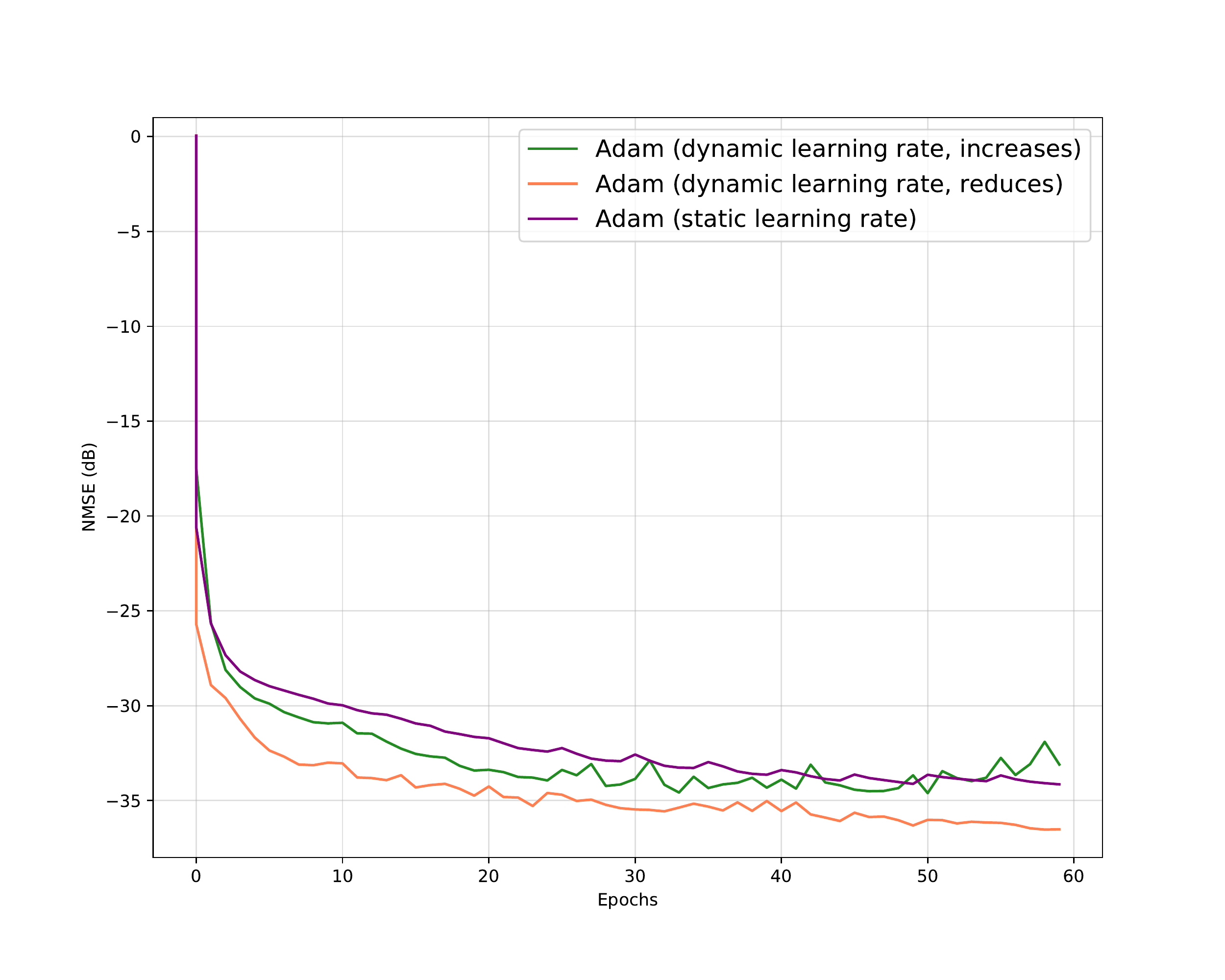}
    }
    \hspace{0.02\textwidth}
    \subfloat[The dependency of NMSE value (dB) on the training time (s)]{
        \includegraphics[width=0.47\textwidth]{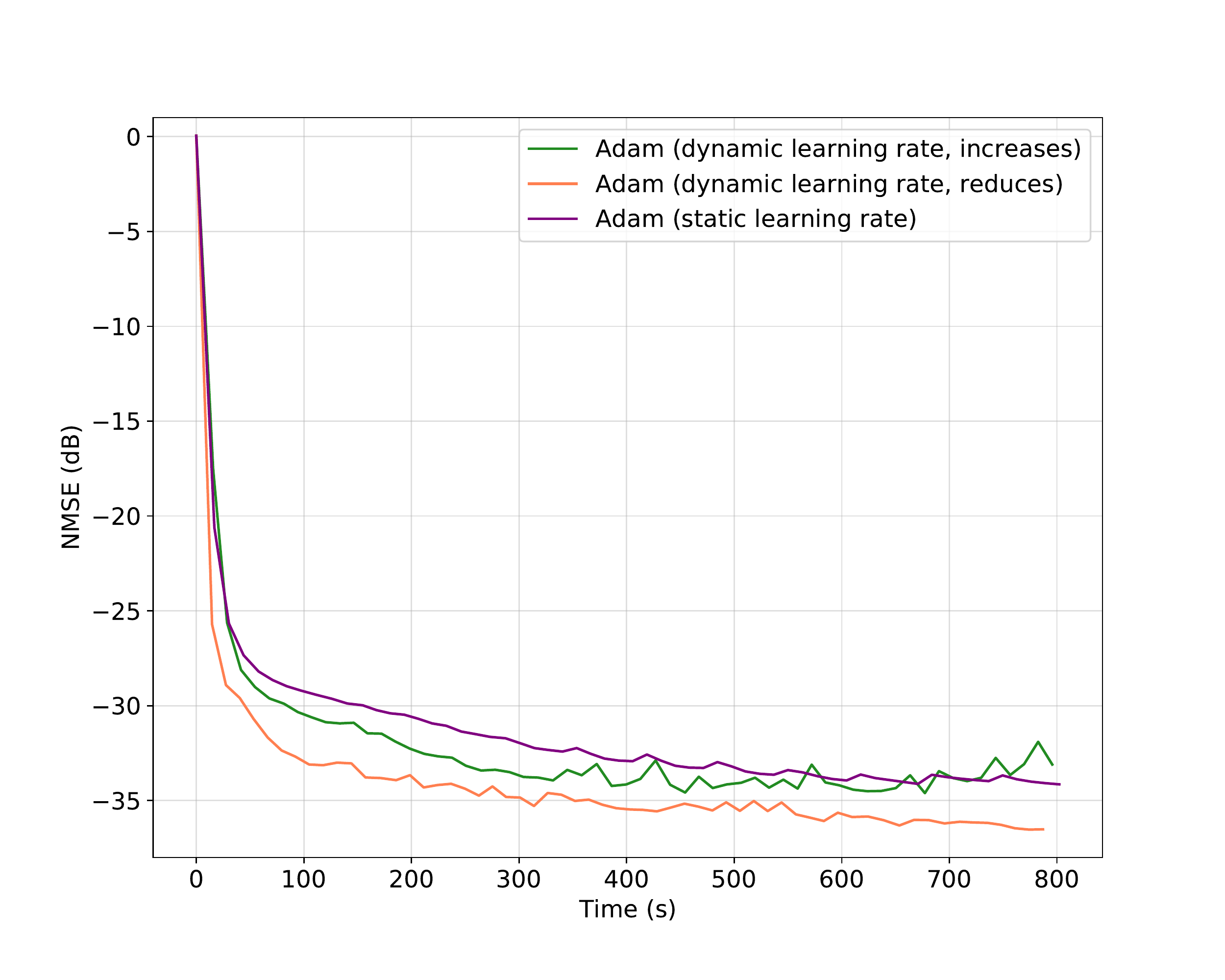}
    }
    \caption{The effect of learning rate changing on training the DPD model. When decreasing, it leads to the significantly more deep convergence of Adam}
    \label{fig:dynamic_lr_static_batch}
\end{figure}

The first group of schedulers consists of the simple monotonically-decaying ones, and here we found the most suitable option for the growth rate and configuration of scheduler. We tested both the increasing and decreasing changing of learning rate, and obtain the following results:

\begin{enumerate}
    \item The most proper formula to change the learning rate is a linear one (bounded below): $\text{lr}~=~\max\{10^{-2} - 10^{-4} \cdot \text{epoch},~6 \cdot 10^{-3}\}$. The effect of scheduling is nearly the same independent of law it follows, if the lower and upper bounds for learning rate are that ($6 \cdot 10^{-3}$ and $10^{-2}$, correspondingly).
    
    \item The learning rate increasing strategy leads to a destabilization of convergence (as expected).
    
    \item The reference score, achieved by the Adam without learning rate scheduling, is $\text{NMSE} = -34.29$~dB. For the properly scheduled learning rate this score becomes $\text{NMSE} = \ok{40}{-36.56}$~dB. Thus, the improvement is 7\%. The convergence curves of some tested options are presented in Fig.~\ref{fig:dynamic_lr_static_batch}
\end{enumerate}

The results summarized above are related to only the monotone schedulers. On the other hand, it becomes more popular to use some cyclic strategies in some cases to improve the properties of obtained solution. For instance, in \cite{Izmailov} it was proposed the scheduler called SWA (Stochastic Weight Averaging), changing learning rate in cyclic manner and averaging the points after the smaller steps. In many practical deep learning problems it allows to achieve around 1\% improvement in accuracy even without a precise tuning of scheduler. The behaviour of cyclic learning rate scheduling without averaging is also considered in \cite{Smith}. We tested both SWA and simple cycling strategies to compare their efficiency on our problem.

\begin{figure}[H]
    \centering
    \subfloat[The behaviour of SWA near the local minimum in many deep learning problems cases. Picture is similar to that from \cite{Izmailov}]{
        \includegraphics[width=0.38\textwidth]{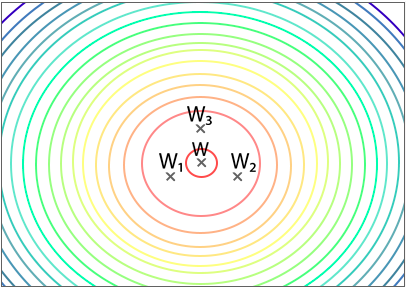}
        \label{fig:swa-expected}
    }
    \hspace{0.02\textwidth}
    \subfloat[The probable behaviour of SWA in case of more complicated structured local minimum.]{
        \includegraphics[width=0.40\textwidth, cfbox=black 0.1pt]{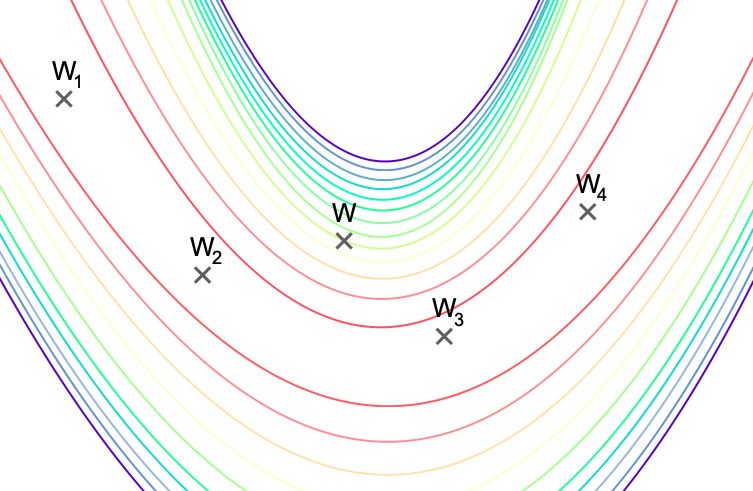}
        \label{fig:swa-in-fact}
    }
    \caption{The two probable cases of SWA \cite{Izmailov} operation. We denote the result of averaging intermediate points $w_1$, $w_2$, $w_3$ by $w$. The left figure demonstrates the simple task for SWA, the right one maybe explained the losing of SWA we observed.}
\end{figure}

\begin{enumerate}
    \setcounter{enumi}{3}
    
    \item The SWA scheduler( with one epoch of fixed learning rate, period of learning rate decreasing equal to one epoch, and same lower and upper bounds on learning rate as in previous experiment) worsens convergence in comparison with reference score obtained by Adam without any scheduling.
    
    \item The same cyclic scheduling without averaging the points, in opposite, allows to obtain significant improvement.
    
    \item For the Adam optimization method, the improvement in depth given by cyclic learning rate is 5\%. For the Adamax, it is 7\%. The results of experiment are summarized in Table~\ref{tab:summary_cyclic_lr}.

\begin{table}[H]
    \centering
    \begin{tabular}{l|c|c}
    NMSE (dB) & Adam   & Adamax \\ \hline
    fixed  & -33.30 & -29.66 \\
    cyclic & \ok{40}{-35.00} & \ok{40}{-31.86}
    \end{tabular}
    
    \caption{Summary of the result of training the DPD model with fixed and with cyclically changing learning rate}
    \label{tab:summary_cyclic_lr}
\end{table}
    
    \item We suppose that cyclic scheduling for our problem indeed behaves like it described in \cite{Garipov} (one of the main prerequisite papers for SWA). But averaging itself fails and it is not asymptotically equivalent to FGE algorithm from \cite{Garipov} in fact. It may be because the structure of local minima in considered problem is not so well as in common deep learning tasks, i.e. instead of situation like in Fig.~\ref{fig:swa-expected} we have something like Fig.~\ref{fig:swa-in-fact}, so we see the loss in score after averaging.
\end{enumerate}

Finally, it is reasonable to assess the learning rate scheduling techniques we considered in compound with the previously described batch size changing, and check if we can combine them to obtain the better scores. We considered to test only a paired technique of increasing the batch size and monotone decreasing of learning rate, since these two are efficient and simple enough to be combined without additional synchronizations.

\begin{figure}[H]
    \centering
    \subfloat[The dependency of NMSE value (dB) on the number of training epochs]{
        \includegraphics[width=0.47\textwidth]{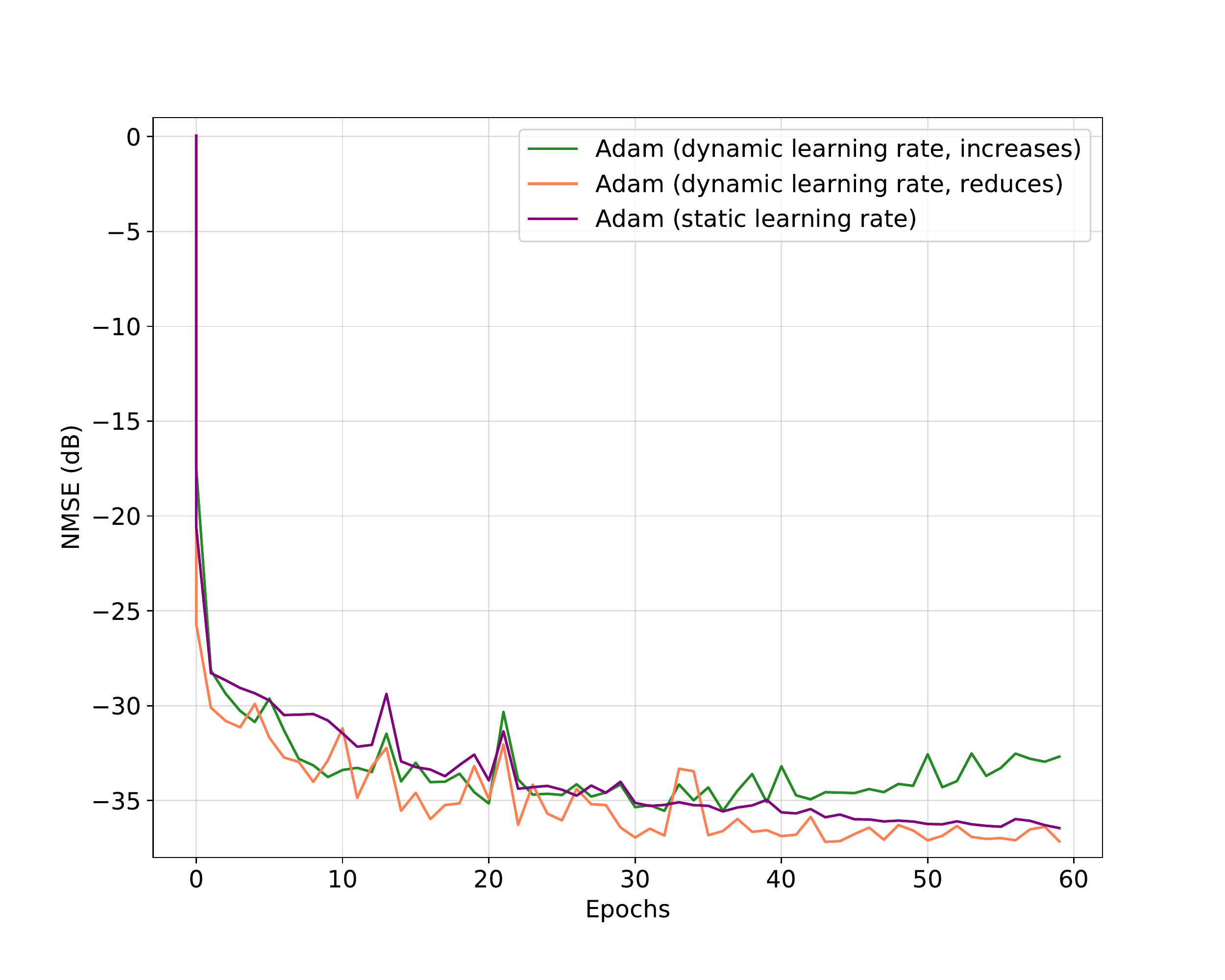}
    }
    \hspace{0.02\textwidth}
    \subfloat[The dependency of NMSE value (dB) on the training time (s)]{
        \includegraphics[width=0.47\textwidth]{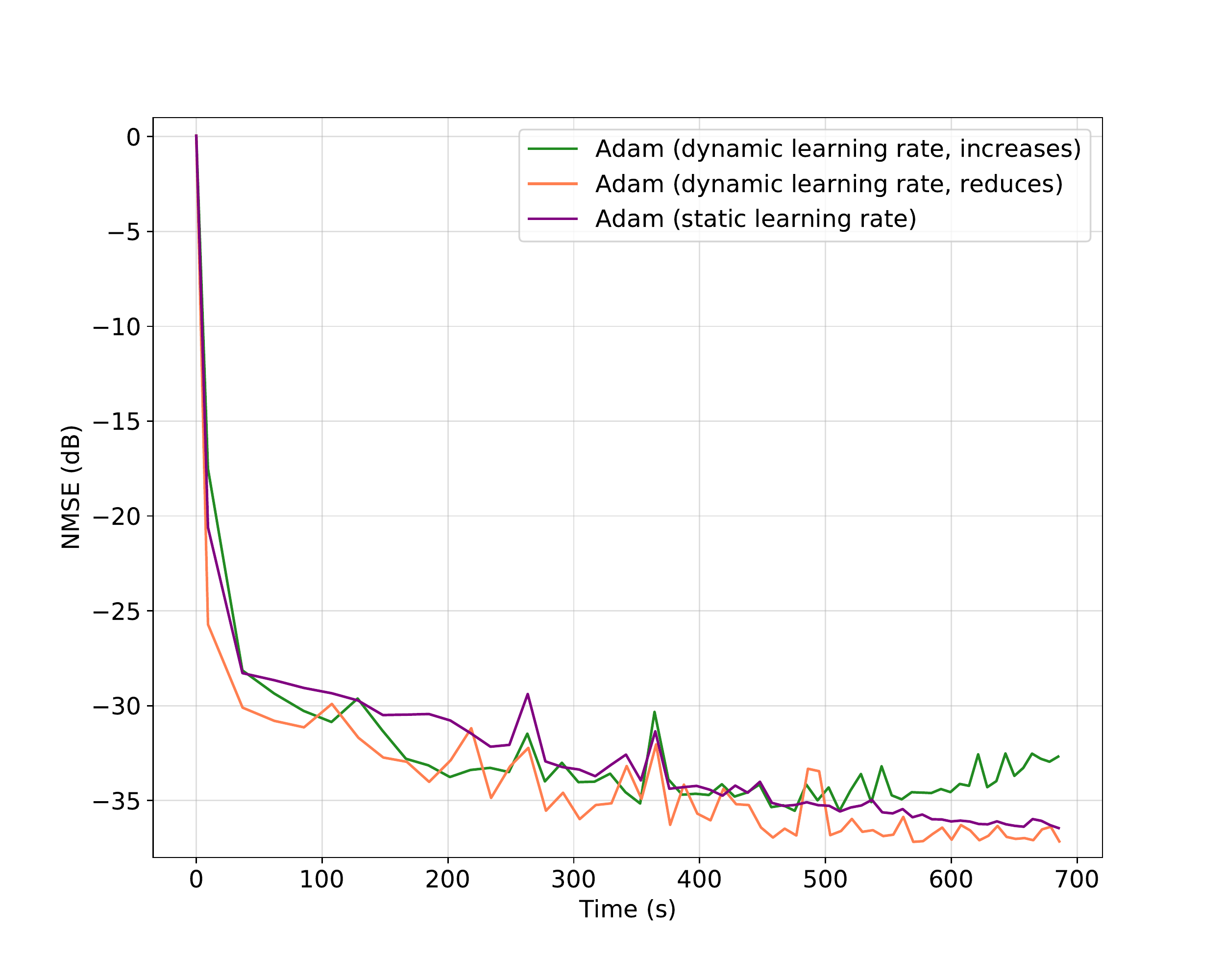}
    }
    \caption{The effect of both learning rate changing and batch size increasing on training the DPD model. When learning rate decreasing, it leads to a little more deep convergence of Adam, while if increasing, it leads to a significant destabilization}
    \label{fig:dynamic_lr_dynamic_batch}
\end{figure}

\begin{enumerate}
    \setcounter{enumi}{7}
    
    \item Increasing batch size with decreasing learning rate is efficient, but the effect of learning rate here is weaker then without batch size changing: it is only a 0.5\%, see Fig.~\ref{fig:dynamic_lr_dynamic_batch}. This is near the scale of score oscillations from epoch to epoch, so this improvement may be unstable.
    
    \item For this compound regime the resulting score is $\text{NMSE} = \indif{40}{-35.46}$ dB. The reference value for this experiment is when batch size increases but learning rate is fixed, it is $\text{NMSE} = -35.27$ dB.
    
    \item Increasing learning rate still leads to a divergence, but this effect become more pronounced when batch size increases.
\end{enumerate}

These results are somewhat unexpected. One can think that if batch size becomes bigger, the stability of method becomes better, and at least with decreasing learning rate (that is very careful heuristic) it should have manifested itself in the improvement of scores, but it did not. Anyway, we can say, that combining of batch size changing with learning rate scheduling is not very efficient, so in any particular case it will be necessary to make choice between big improvement in time and good improvement in depth. It mentioned, that working time is minor metric for signal processing tasks. Nevertheless, it would be interesting to know that in online regime changing batch size is just more efficient with respect to the depth of convergence than both of learning rate scheduling policies.

\subparagraph{Regularization}
The last of non-optimizational modifications considered in the course of our study is a regularization. It is very common to add the regularization term in the problem to control the properties of the solution, but the regularization may also be a good heuristic to improve the convergence, due to its impact on (local) conditionality of the problem. This practical property of regularization seems to be underexamined and is bypassed in many new research papers on application of optimization methods to learning problems.

\begin{figure}[H]
    \centering
    \subfloat[The dependency of NMSE value (dB) on the number of training epochs for Adam]{
        \includegraphics[width=0.47\textwidth]{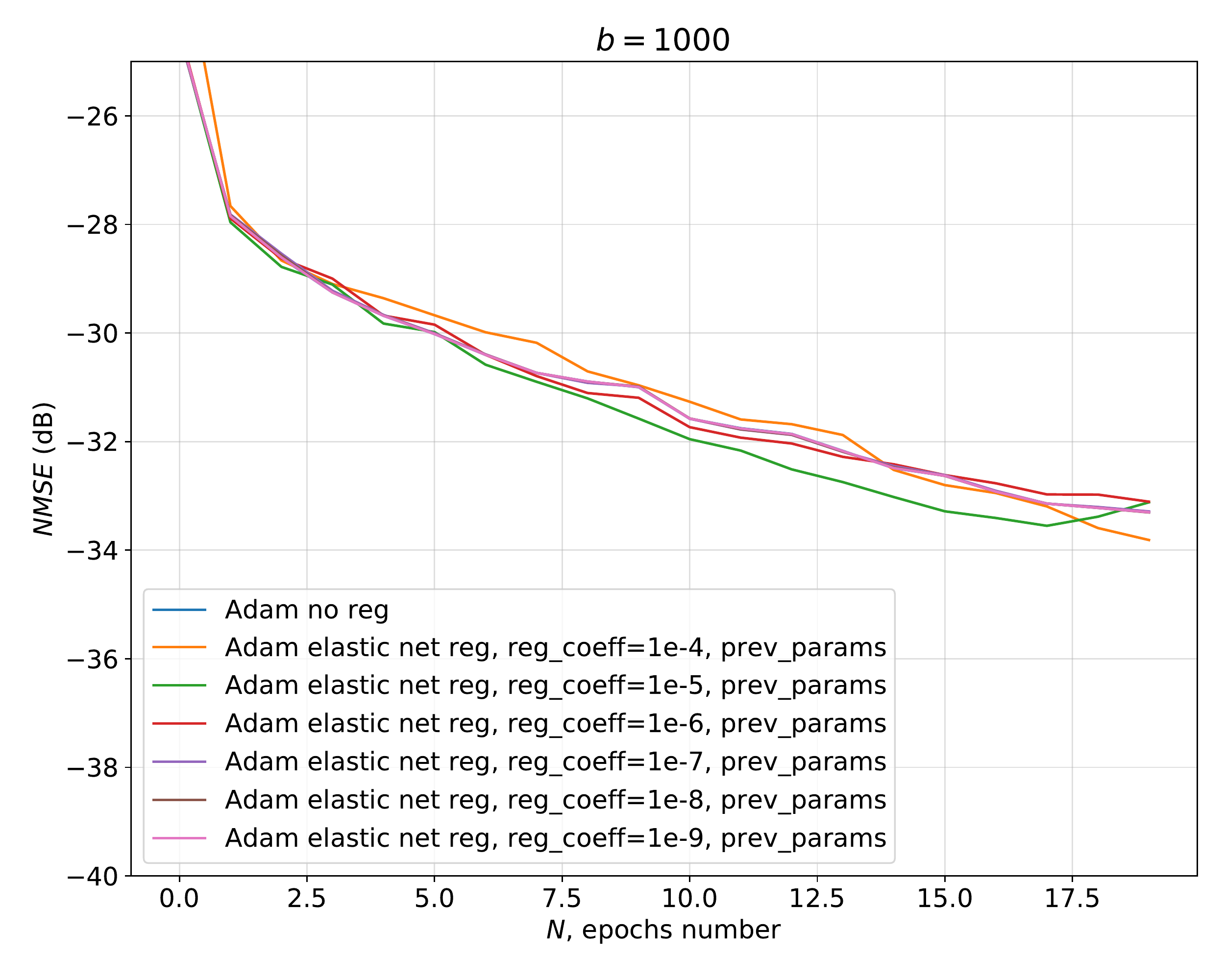}
    }
    \hspace{0.02\textwidth}
    \subfloat[The dependency of NMSE value (dB) on the number of training epochs for RAdam]{
        \includegraphics[width=0.47\textwidth]{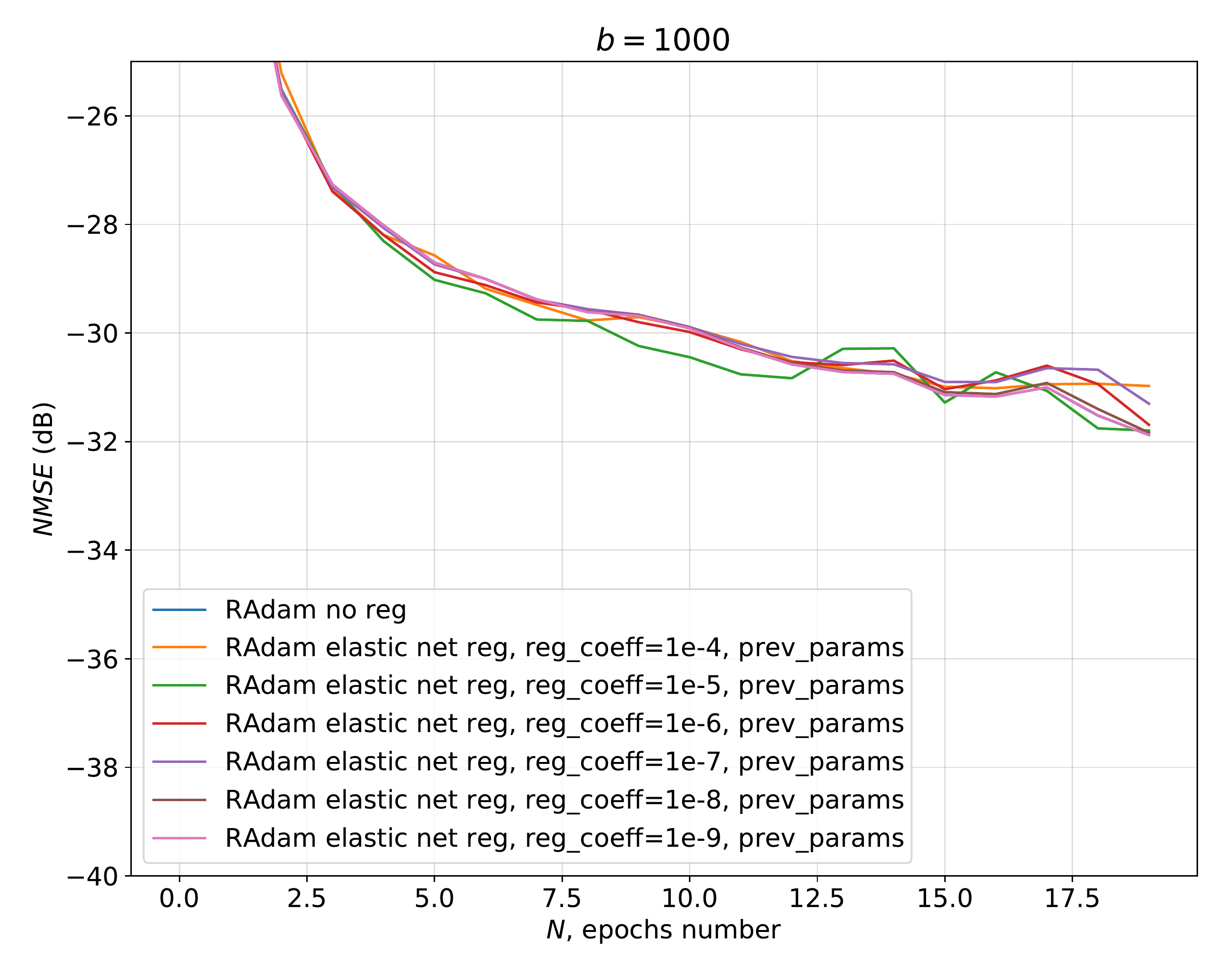}
    }
    \caption{The effect of prox-centered regularization on training the DPD model}
    \label{fig:nmse_reg_prev}
\end{figure}

In out experiment, we considered only some widely-known regularization terms, namely $\ell_1$~(lasso), $\ell_2$~(ridge) and $\ell_1 + \ell_2$~(elastic net). To check their impact on convergence, we tested the efficiency of different optimization methods, listed in Table~\ref{tab:methods}, on the corresponding regularized problems. And there are the results it gave:

\begin{enumerate}
    \item Unlike the most of results presented in previous sections, the impact of regularization is significantly dependent on the chosen optimization method. So, Adam and Adamax are sensitive to it, while RAdam \cite{Liu} does not totally (see Table~\ref{tab:reg}).
    
    \item On the other hand, one can conclude from the left part of Table~\ref{tab:reg} that the standard regularization (one like $\lambda \|\theta\|^2_2$, we name it zero-centered; we set $\lambda = 10^{-4}$) does not affect the convergence. We present the minimum-aggregated values of NMSE, so we hide some probable destabilization effects here~--- they can appear in RAdam, but not in Adam or Adamax (see Fig.~\ref{fig:nmse_reg_prev}).
\end{enumerate}

The point is that it more efficient here is to use more complicated regularization with moving center (one like $\lambda \|\theta - \overline{\theta}\|^2_2$, we name it prox-centered). We tested the version with a period of updates equal to one epoch and $\overline{\theta}$ set to point obtained after the last iteration of method. This configuration resembles the scheme of inexact proximal gradient method a little.

\begin{enumerate}
    \setcounter{enumi}{2}
    
    \item So, this approach is efficient and gives 2\% improvement to the Adam optimization method. The results are summarized in Table~\ref{tab:reg}.
    
    \item It is unusual that the best result is obtained with $\ell_1$-regularization (and also with elastic net, but the effect is smaller, so we can conclude that $\ell_2$ part of it only suppresses the effect of $\ell_1$ one).
\end{enumerate}

\begin{table}[H]
    \centering
    \begin{tabular}{l|c|c|c|c|c|c|c}
    \multirow{2}{*}{method} & \multirow{2}{*}{reference} & \multicolumn{3}{c|}{min NMSE, zero-centered} & \multicolumn{3}{c}{min NMSE, prox-centered} \\ \cline{3-8}
    & & ${\ell_1}$ & ${\ell_2}$ & ${\ell_1+\ell_2}$ & ${\ell_1}$ & ${\ell_2}$ & ${\ell_1+\ell_2}$ \\ \hline
    Adam  & -33.30 & \indif{40}{-33.30} & \indif{40}{-33.30} & \indif{40}{-33.30} & \ok{40}{-33.99} & \indif{40}{-33.31} & \ok{40}{-33.81} \\
    RAdam & -31.88 & \indif{40}{-31.88} & \indif{40}{-31.88} & \indif{40}{-31.88} & \indif{40}{-31.88} & \indif{40}{-31.88} & \indif{40}{-31.88}
    \end{tabular}
    
    \caption{Summary of the result of training the DPD model with different regularization terms}
    \label{tab:reg}
\end{table}

\paragraph{Quasi-online learning framework}

\begin{figure}[H]
\centering
\includegraphics[width=0.6\linewidth]{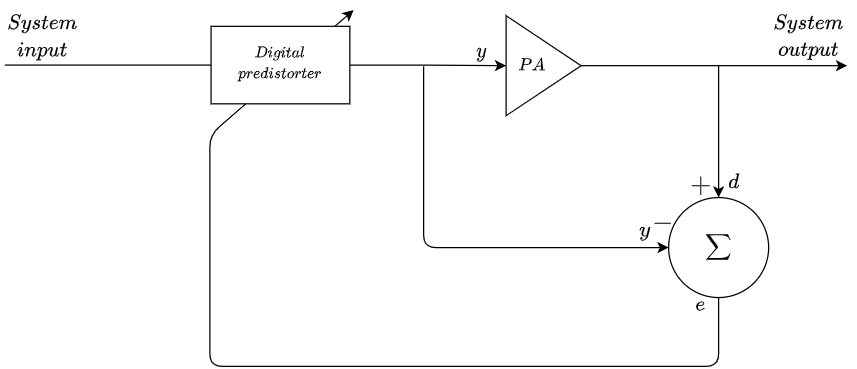}
\caption{Schematic diagram of the operation and learning of the DPD device during signal processing}
\label{fig:dpd}
\end{figure}

It was said before that the standard learning framework is not suitable for learning the DPD models due to some disagreements with the practical devices behaviour. The point is that if we use some optimization methods, for example, Adam and Adamax, to fit the parameters of the real technical prototype, and then test it on some real life signal, we will get that Adamax (significantly) outperforms Adam, while our experiments in previous section predict the opposite. The reason is that the operation of technical prototype is different from the simulating scheme used before, due to the specifics of signal processing task. In fact, the DPD model retrains periodically to be able to work with the new incoming signal. Unlike machine learning, the data in this case is updated with a big frequency (80 MHz or more) and, moreover, change its properties, for example, due to the change of the interlocutor on the telephone line. A very simplified circuit diagram of DPD operation is presented in Fig. \ref{fig:dpd}, and represents the looped process. So, if our results disagree with the experiment, we firstly should take into account this cyclicality as the most notable difference with the standard machine learning framework.
Further, we propose another framework to train/test DPD models, that includes the retraining procedures and more correct evaluation, and test some of the techniques proposed above in new, more close to real, conditions.

\subparagraph{New framework and comparison of algorithms}
As it was noted in the preamble of section, the provided comparison of methods and testing any practices in the previously postulated conditions have little in common with the behaviour of these objects observed in reality. At least, it is difficult to find real analogies for training and test datasets, and yet the NMSE value we use is a function not only of the applied numerical method, but of these two parameters too. Of course, we are not able to avoid this functional dependency. We propose to replace these two variables with the most natural one~--- the whole given signal segment. 

\begin{figure}[H]
    \centering
    \includegraphics[width=0.7\linewidth]{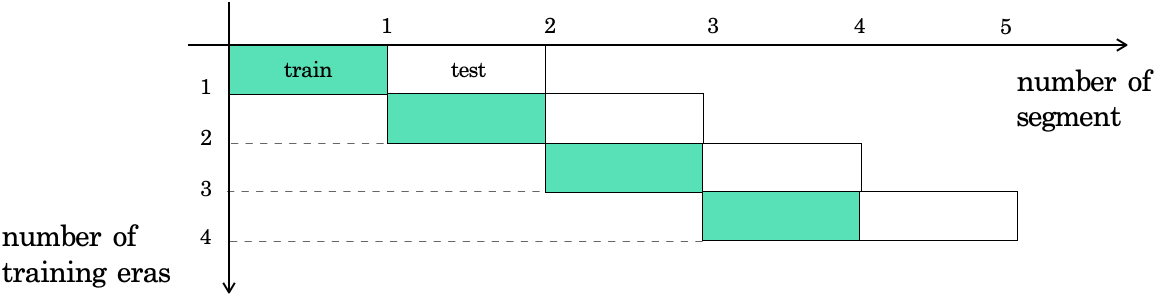}
    \caption{Description of an online learning framework: a diagram of data segments used at different moments in signal processing by the model}
    \label{fig:online_description}
\end{figure}

Let us turn to the operation of DPD device. Following the scheme from Fig.~\ref{fig:dpd}: DPD collects information about some small segment of signal, then it trains on it and is predicting distortion for some new small segment, while collecting information about it, then it repeats this procedure over and over again. The only thing we cannot wave away here is ``over and over again'', because we are limited with the given dataset, but we can exactly simulate all the other steps. Namely, let us split our dataset into some number of segments with the same length and train and test model on the two consequent segments correspondingly at the every training \textit{era}. This is shown in Fig.~\ref{fig:online_description}.

There are two options to assess model in this framework. The most obvious one is the value of NMSE from the last era, i.e. evaluated on the last segment of dataset. But this metric has two flaws: 1) it is evaluated only on the last segment of data, that may lead to an approximation issues (like overestimation of error), if this segment is small; 2) it does not represent the performance of the method the whole signal, and even if all segments except the last one were predicted absolutely incorrect, score can be good enough. To deal with these problems, we also use the mean-aggregated version of NMSE:
\begin{equation*}
   \text{mean NMSE}(T, y, \overline{y}):= \frac{1}{T} \sum_{i=1}^T NMSE(y_i, \overline{y}_i) \text{ dB}.
\end{equation*}
It is in fact a very proper choice for assessment the performance of online optimization methods, because it is usually used in their theoretical analysis \cite{Hazan} (and is connected with so called regret metric). 

\begin{figure}[H]
    \centering
    \subfloat[The dependency of NMSE value (dB) on the working time (s)]{
        \includegraphics[width=0.47\textwidth]{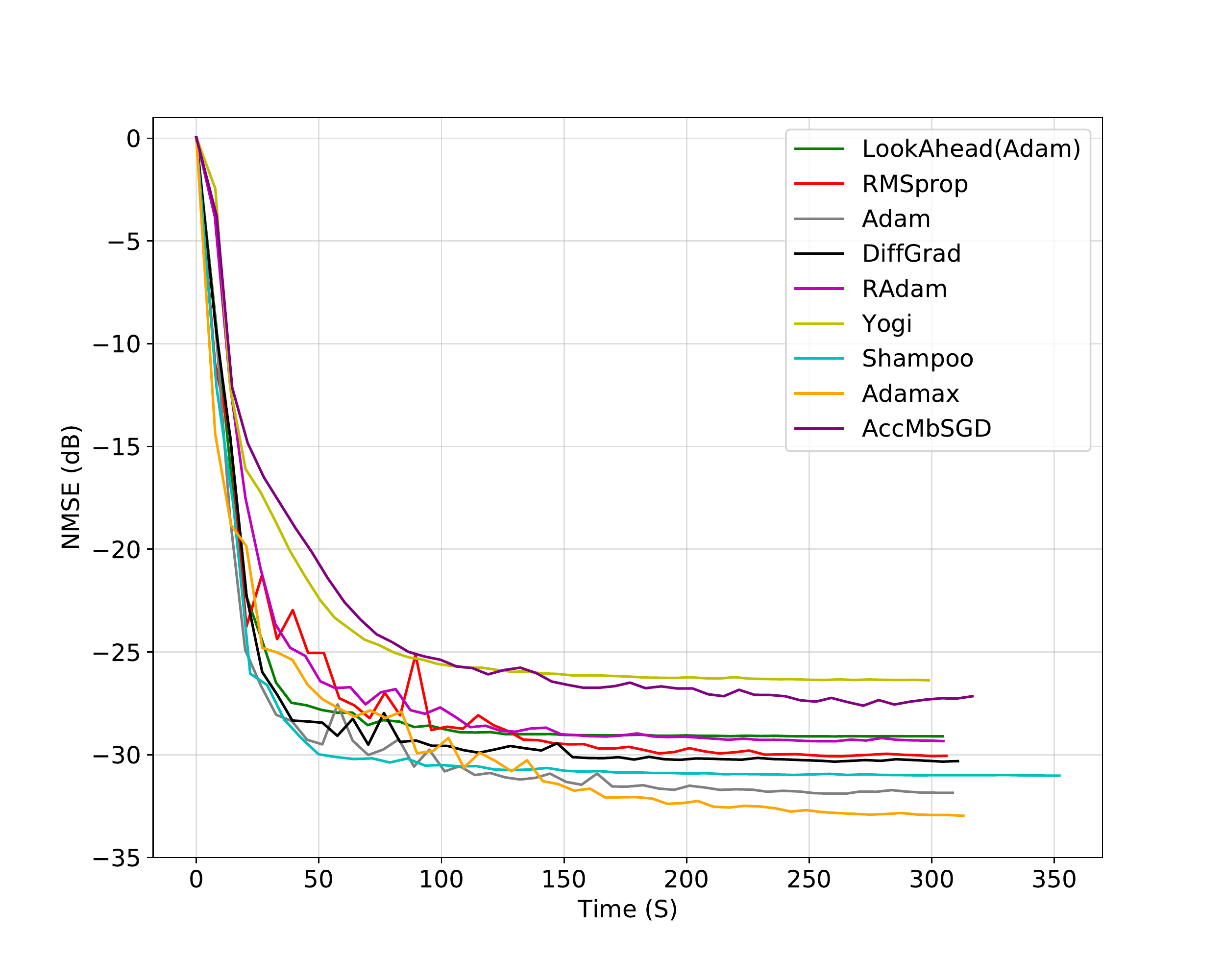}
        \label{fig:o_methods}
    }
    \hspace{0.02\textwidth}
    \subfloat[The dependency of mean time NMSE value (dB) on the working time (s)]{
        \includegraphics[width=0.47\textwidth]{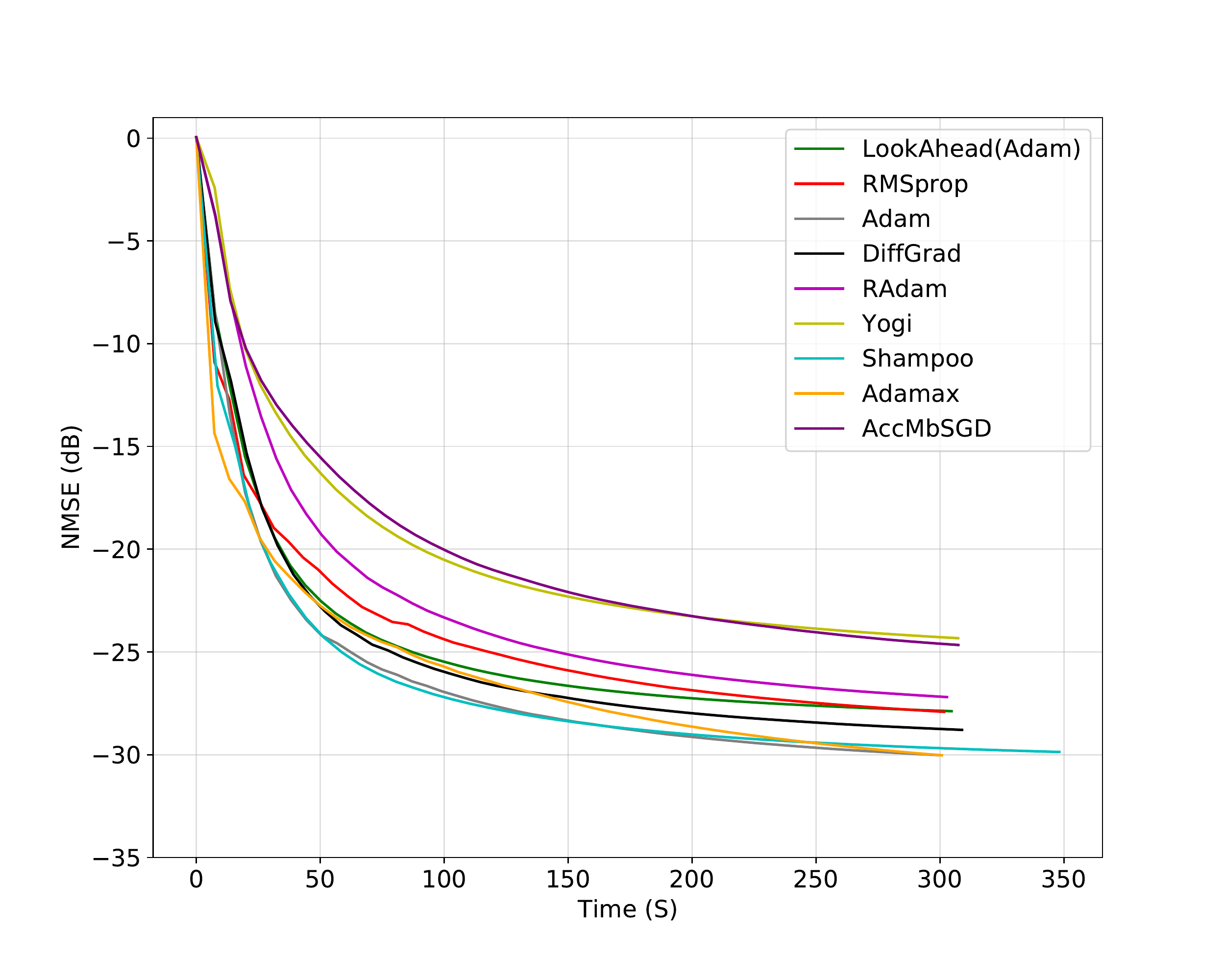}
        \label{fig:o_methods_mean}
    }
    \caption{The convergence curves of different optimization methods in quasi-online regime with evaluation after every era}
\end{figure}

So, let us firstly reassess the optimization methods in new setting without any additional modifications. In Table~\ref{tab:methods} we present the results for some of the methods with corresponding values of NMSE and mean NMSE. And here we see the key difference of new approach with the previous one: Adamax algorithm now is the best performing method and by a large margin~--- it is 3\% better than Adam in NMSE, and also have better rate of convergence at the late eras (see Fig.~\ref{fig:o_methods}). The latter is the most important property of the optimization method in our task, because it determines the possible depth of solution given by this method. We see, that despite the starting lag of Adamax it managed to overtake methods like Shampoo and Adam, that seems to be stick in local minima or just have decline in convergence rate.

\begin{table}[H]
    \centering
    \begin{tabular}{l|c|c}
    Method   & NMSE (dB) & mean NMSE (dB) \\ \hline
    Adamax   & \ok{40}{-33.04} & \ok{40}{-30.21} \\
    Adam     & -31.94   & -30.17 \\
    Shampoo  & -31.06   & -29.99 \\
    DiffGrad & -30.41   & -28.76 \\
    RMSprop  & -30.18   & -27.96 \\
    LaAdam   & -29.12   & -27.90 \\
    RAdam    & -29.35   & -27.24 \\
    AccMbSGD & -27.18   & -24.65 \\
    Yogi     & -26.41   & -24.29
    \end{tabular}
    
    \caption{Summary of the result of training the DPD model with different optimization methods in new quasi-online setting}
    \label{tab:methods}
\end{table}

For the other methods we see mostly the same results as that for standard learning framework. Shampoo is still on of the best starter methods, and AccMbSGD demonstrates not well convergence rate. The Fig.~\ref{fig:o_methods_mean} shows the mean-aggregated version of the same plots, and we can make the similar conclusions from it, except that leadership of Adamax here is not so pronounced, but it is just an effect of averaging. 

\subparagraph{Dynamic batch size}
Now, we should revisit the practices we proposed for the standard learning framework in the view of new one. Let us start with the batch size changing. Further we consider only the most efficient policy we identified above (linear increasing of batch size), but at the same time, we test two modifications of this scheme, that are given by the more complicated structure of testing framework itself.

Namely, the first option is to change batch size linearly in dependence on current era (number of learning segment), and the second option is to bind batch size with the number of epoch within every of the segments every era. The choice of policy here is in fact the trade-off between depth of convergence in ``expected'' score and depth of convergence for an every segment correspondingly.

We tested the both policies, and the result is that era changing (first option) and epoch changing (second option) are competitive depending of what the metric we use. In the standard NMSE, era changing performs better, but after averaging epoch changing becomes a better option. It means, that it is more efficient to achieve a little improvement in convergence at the every era (for the every small piece of signal) instead of improving the outer loop convergence. It is quite unexpected result. 

\begin{figure}[H]
    \centering
    \subfloat[The dependency of NMSE value (dB) on the working time (s)]{
        \includegraphics[width=0.47\textwidth]{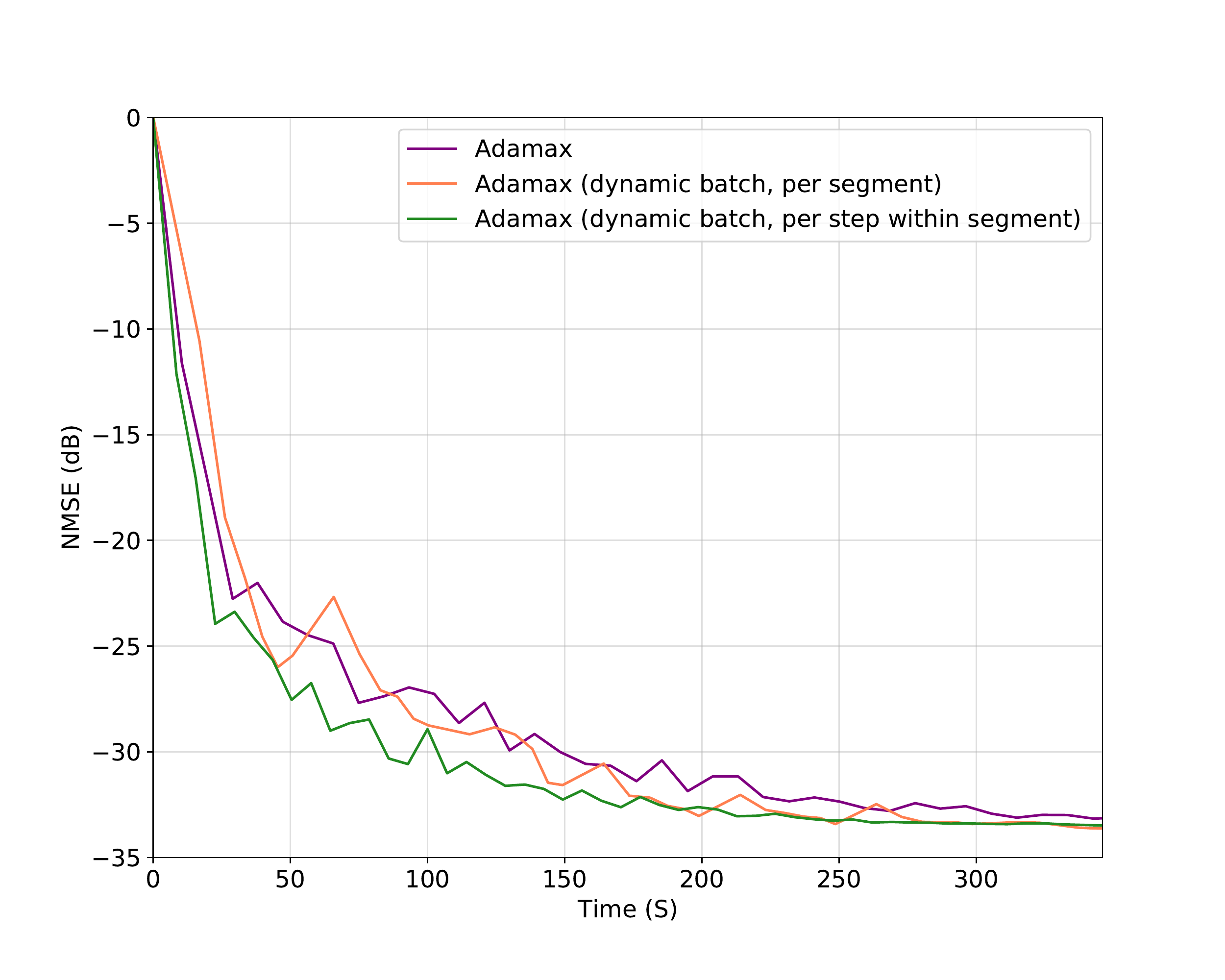}
    }
    \hspace{0.02\textwidth}
    \subfloat[The dependency of mean time NMSE value (dB) on the working time (s)]{
        \includegraphics[width=0.47\textwidth]{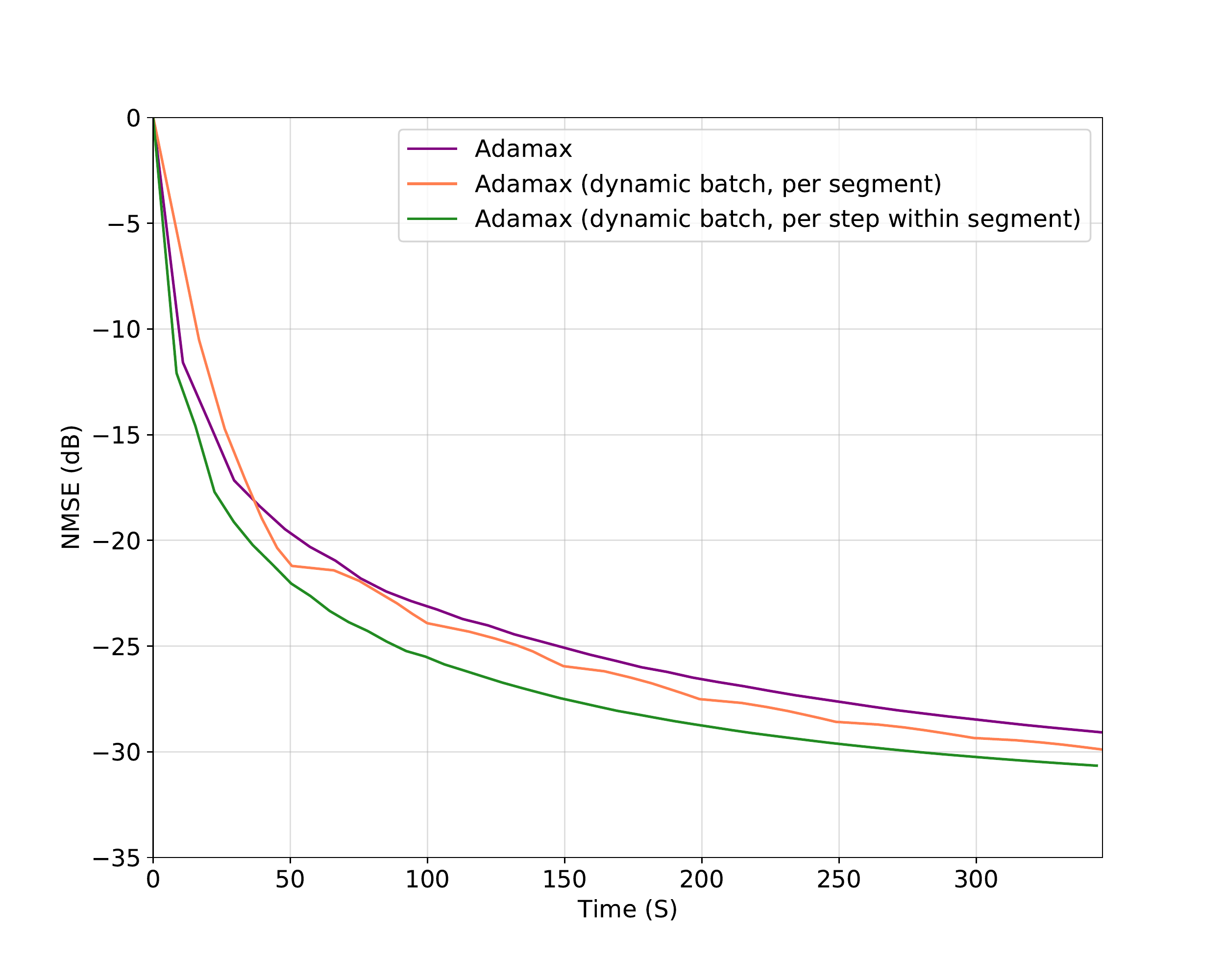}
    }
    \caption{The effect of batch size changing on training the DPD model}
    \label{fig:convergence_o_dynamic_batch_size}
\end{figure}

The results of this experiment are summarized in Table~\ref{tab:summary_o_dynamic_batch}. The use of the epoch changing batch size allows to obtain 6\% improvement in mean NMSE metric for the considered setting. Fig.~\ref{fig:convergence_o_dynamic_batch_size} shows the convergence of Adamax with the different batch size changing policies in more details for two metrics.

\begin{table}[H]
    \centering
    \begin{tabular}{l|c|c}
                  & NMSE (dB) & mean NMSE (dB) \\ \hline
    fixed        & -33.15    & -28.96  \\
    era changing   & \ok{40}{-33.63} & \ok{20}{-29.85} \\
    epoch changing & \ok{20}{-33.49} & \ok{40}{-30.66}
    \end{tabular}
    
    \caption{Summary of the result of training the DPD model with fixed and with changing batch size in quasi-online regime}
    \label{tab:summary_o_dynamic_batch}
\end{table}

\subparagraph{Regularization}
In the case of regularization, we do not need to take into account the difference between eras ans epochs, so to test the performance of this modifications we just carried out one more experiment of learning DPD model with Adamax optimization method (it is the best one in online setting) with $\ell_1$-regularization (that was the best and the most stable option in standard learning framework experiments). The result approves the positive impact of $\ell_1$-regularization: comparing to the reference value of $\text{mean NMSE} = -30.72$ dB, given by the model trained without any regularization, $\ell_1$-regularization allows to achieve $\text{mean NMSE} = \indif{40}{-30.81}$ dB, that is 0.3\% better. Unfortunately, this is near the scale of oscillations in convergence on the late iterations of the method, so we conclude that regularization gives no improvement in online regime (but also does not impair the convergence depth). 

We adhere to the explanation of this that the oscillations in online regime are itself really perceptible, and the effect of our modification is just suppressed by multiple switching between different data segments. Another reason is that after the last change of segment, we again start to seek new optimal point (corresponding to this segment) from some arbitrary point, and due to the small value of regularization coefficient, we have the situation similar to that from standard framework experiment, but within the small subset of data. In other words, regularization impacts convergence only on the last segment, and impact is small because of smallness of this segment.

\subparagraph{Cyclic learning rate}
At last, we reproduced the experiments with cyclic learning rate scheduling in the new online framework. Now we applied this modification to the new best method, Adamax. We tested both options for the period of cycling: once in epoch or once in era, and chose the epoch~--- there is no qualitative difference between the results in two these configurations, but increasing the period leads to a decrease of this positive effect itself. So, we used learning rate cycling scheduling with the period equal to one epoch together with Adamax, for which the improvement was 1\% in mean NMSE, and with Adam, with the similar resulting 1\% improvement in mean NMSE. The decrease in effect in comparison with that one for standard learning framework is explained by the same reasons as presented in previous paragraph.

\paragraph{Conclusion}

In the introduction, we posed two central questions of this research: which optimization method is the most suitable for training the DPD model, and what the framework we need to use to simulate the DPD operation in more close to real conditions and to agree our simulation with the results observed in practice. In the course of our study, we carried out a big series of clarifying and detailing experiments to come to the two desired answers. In the first section, dedicated to the standard learning framework, we summarize those that relate to the first question, presenting the comparison of many modern optimization methods and revealing the gist of proposed modifications designed to improve their convergence. In the second section, we propose the new quasi-online learning framework and demonstrate that the ranking of the methods in our new experiments is mostly in line with that observed in close to real conditions, and also reassess some of proposed modifications to approve their stability.

To express the achieved improvements in numbers: in standard learning framework, the depth of convergence (value of NMSE metric, that is logarithmic scaled itself) was improved in 7\% (thanks to steadily decreasing learning rate changing for Adam or cyclic scheduling for Adamax), while the maximum achieved improvement in online regime is 6\% in mean NMSE (and given by the increasing of batch size at every epoch). In terms of working time, we have achieved a twofold improvement (with changing of batch size), that preserves 3\% and 6\% improvement in NMSE and mean NMSE for the standard and online regime, correspondingly.

Interim results arising in the course of the study confirm that the considered DPD model is really specific and requires a separate study of the effectiveness of optimization methods as applied to its training: many approaches, common for example in the field of training neural networks, work here unexpectedly or do not work at all. The author's team hopes that the proposed results will be useful for other researchers of DPD models, and that the questions asked will stimulate a deeper study of problems in the vicinity of this work.

\paragraph{Acknowledgments}

The work on this research was supported and organized in cooperation with the Educational Сenter ``Sirius'' and Educational Foundation ``Talent and Success'', Russia, during the Project science and technology program
``Big Challenges''. The  assistance provided by the mentor of program, Anton Gusev, was greatly appreciated. The topic of research is a part of joint project of Moscow Institute of Physics and Technology and Huawei Russian Research Institute. We would like to thank Andrey Vorobyev, Eugeniy Yanitskiy and Anna Nikolaeva for openness to cooperation and interest in this study. We wish to extend out special thanks to prof. Alexander Gasnikov for all the helpful discussions and his unwavering favor for this project. We also thank Liliana Baislanova and Mansur Zainullin for their participation in setting up some numerical experiments.


\begin{thebibliography}{99}
\providecommand{\url}[1]{\texttt{#1}}
\providecommand{\urlprefix}{URL }
\providecommand{\doi}[1]{https://doi.org/#1}

\bibitem[Haykin, 2008]{Haykin} {\it Haykin~S.\,S.}
    Adaptive filter theory.~--- Pearson Education India, 2008.
    
\bibitem[Ghannouchi, Hammi, Helaoui, 2015]{Ghannouchi} {\it Ghannouchi~F.\,M., Hammi~O., Helaoui~M.}
    Behavioral modeling and predistortion of wideband wireless transmitters.~--- John Wiley \& Sons, 2015.

\bibitem[Schreurs et al., 2008]{Schreurs} {\it Schreurs~D. et al.}
    RF power amplifier behavioral modeling.~--- New York, USA: Cambridge university press, 2008.

\bibitem[Pasechnyuk et al., 2021]{Pasechnyuk} {\it Pasechnyuk~D. et al.}
    Non-convex optimization in digital pre-distortion of the signal~// arXiv preprint arXiv:2103.10552.~--- 2021.
    
\bibitem[Kingma, Ba, 2014]{Kingma} {\it Kingma~D.\,P., Ba~J.}
    Adam: A method for stochastic optimization~// arXiv preprint arXiv:1412.6980.~--- 2014.
    
\bibitem[Woodworth, Srebro, 2021]{Woodworth} {\it Woodworth~B., Srebro~N.}
    An Even More Optimal Stochastic Optimization Algorithm: Minibatching and Interpolation Learning~// arXiv preprint arXiv:2106.02720.~--- 2021.

\bibitem[Gupta, Koren, Singer, 2014]{Gupta} {\it Gupta~V., Koren~T., Singer~Y.}
    Shampoo: Preconditioned stochastic tensor optimization~// International Conference on Machine Learning.~--- PMLR, 2018.~--- p. 1842-1850.

\bibitem[Liu et al., 2019]{Liu} {\it Liu~L. et al.}
    On the variance of the adaptive learning rate and beyond~// arXiv preprint arXiv:1908.03265.~--- 2019.

\bibitem[Zhang et al., 2019]{Zhang} {\it Zhang~M.\,R. et al.}
    Lookahead optimizer: k steps forward, 1 step back~// arXiv preprint arXiv:1907.08610.~--- 2019.
    
\bibitem[Izmailov et al., 2008]{Izmailov} {\it Izmailov~P. et al.}
    Averaging weights leads to wider optima and better generalization~// arXiv preprint arXiv:1803.05407.~--- 2018.

\bibitem[Garipov et al., 2008]{Garipov} {\it Garipov~T. et al.}
    Loss surfaces, mode connectivity, and fast ensembling of dnns~// Proceedings of the 32nd International Conference on Neural Information Processing Systems.~--- 2018.~--- p. 8803-8812

\bibitem[Smith, 2017]{Smith} {\it Smith~L.\,N.}
    Cyclical learning rates for training neural networks~// 2017 IEEE winter conference on applications of computer vision (WACV).~--- IEEE, 2017.~--- p.~464-472.

\bibitem[Zhao, Xie, Li, 2020]{Zhao} {\it Zhao~S.\,Y., Xie~Y.\,P., Li~W.\,J.}
    Stagewise Enlargement of Batch Size for SGD-based Learning~// arXiv preprint arXiv:2002.11601.~--- 2020.

\bibitem[Devarakonda, Naumov, Garland, 2017]{Devarakonda} {\it Devarakonda~A., Naumov~M., Garland~M.}
    Adabatch: Adaptive batch sizes for training deep neural networks~// arXiv preprint arXiv:1712.02029.~--- 2017.

\bibitem[Xiong et al., 2020]{Xiong} {\it Xiong~R. et al.}
    On layer normalization in the transformer architecture~// International Conference on Machine Learning.~--- PMLR, 2020.~--- p. 10524-10533.

\bibitem[Hazan, Kale, 2014]{Hazan} {\it Hazan~E., Kale~S.}
    Beyond the regret minimization barrier: optimal algorithms for stochastic strongly-convex optimization~// The Journal of Machine Learning Research.~--- 2014.~--- V. 15.~--- No. 1.~--- p. 2489-2512.

\end{thebibliography}
\end{document}